\documentclass[10pt]{article}

\usepackage{amsmath,amsfonts,amsthm,amssymb,latexsym,url}
\usepackage[a4paper,portrait]{geometry}

\newcommand{\EM}{\ensuremath}

\newtheorem{thm}{Theorem}[section]%
\newtheorem{cor}[thm]{Corollary}%
\newtheorem{lem}[thm]{Lemma}%
\newtheorem{xpl}[thm]{Example}%
\newtheorem{rem}[thm]{Remark}%

\newcommand{\dN}{\EM{\mathbb{N}}}

\newcommand{\dP}{\EM{\mathbb{P}}}

\newcommand{\dR}{\EM{\mathbb{R}}}

\newcommand{\dZ}{\EM{\mathbb{Z}}}

\newcommand{\rD}{\EM{\mathrm{D}}}

\newcommand{\rI}{\EM{\mathrm{I}}}

\newcommand{\rL}{\EM{\mathrm{L}}}

\newcommand{\cB}{\EM{\mathcal{B}}}
\newcommand{\cC}{\EM{\mathcal{C}}}
\newcommand{\cD}{\EM{\mathcal{D}}}
\newcommand{\cE}{\EM{\mathcal{E}}}
\newcommand{\cF}{\EM{\mathcal{F}}}
\newcommand{\cG}{\EM{\mathcal{G}}}
\newcommand{\cH}{\EM{\mathcal{H}}}
\newcommand{\cI}{\EM{\mathcal{I}}}

\newcommand{\cK}{\EM{\mathcal{K}}}
\newcommand{\cL}{\EM{\mathcal{L}}}

\newcommand{\cN}{\EM{\mathcal{N}}}

\newcommand{\cP}{\EM{\mathcal{P}}}

\newcommand{\cS}{\EM{\mathcal{S}}}
\newcommand{\cT}{\EM{\mathcal{T}}}

\newcommand{\bE}{\EM{\mathbf{E}}}

\newcommand{\bI}{\EM{\mathbf{I}}}
\newcommand{\bJ}{\EM{\mathbf{J}}}

\newcommand{\bL}{\EM{\mathbf{L}}}

\newcommand{\bR}{\EM{\mathbf{R}}}

\newcommand{\al}{\alpha}
\newcommand{\be}{\beta}

\newcommand{\de}{\delta}
\newcommand{\ga}{\gamma}

\newcommand{\la}{\lambda}

\newcommand{\si}{\sigma}

\newcommand{\te}{\theta}

\newcommand{\veps}{\varepsilon}
\newcommand{\vphi}{\varphi}

\newcommand{\p}[4]{{#3}\!\left#1{#4}\right#2} 

\newcommand{\PENT}[1]{\EM{\lfloor#1\rfloor}} 
\newcommand{\ABS}[1]{\EM{{\left| #1 \right|}}} 
\newcommand{\BRA}[1]{\EM{{\left\{#1\right\}}}} 
\newcommand{\NRM}[1]{\EM{{\left\| #1\right\|}}} 
\newcommand{\PAR}[1]{\EM{{\left(#1\right)}}} 
\newcommand{\pd}{\EM{\partial}} 
\newcommand{\SBRA}[1]{\EM{{\left[#1\right]}}} 

\newcommand{\ENTF}[2]{\entf{#1}^{#2}}
\newcommand{\ENT}[3]{\p(){\ENTF{#1}{#2}}{#3}}
\newcommand{\entf}[1]{\mathbf{Ent}_{#1}}
\newcommand{\ent}[2]{\p(){\entf{#1}}{#2}}
\newcommand{\moyf}[1]{\bE_{#1}}
\newcommand{\moy}[2]{\p(){\moyf{#1}}{#2}}
\newcommand{\varf}[1]{\mathbf{Var}_{#1}}
\newcommand{\var}[2]{\p(){\varf{#1}}{#2}}

\newcommand{\GA}{\boldsymbol{\Gamma}}
\newcommand{\GD}{\GA_{\!\!{\mathbf 2}}}
\newcommand{\GI}{\bL}
\newcommand{\GR}{\nabla}

\newcommand{\SGf}[1]{{\mathbf P}_{\!#1}}
\newcommand{\SG}[2]{\p(){\SGf{\!#1}}{#2}}

\newcommand{\BS}[1]{\EM{\boldsymbol{#1}}}

\renewcommand{\leq}{\leqslant}
\renewcommand{\geq}{\geqslant}
\newcommand{\ds}[1]{\EM{\displaystyle{#1}}}

\title{Binomial-Poisson entropic inequalities and the M/M/$\infty$ queue}

\author{Djalil Chafa\"{\i}}

\date{October 2005. Revised, March 2006.}

\begin{document}

\maketitle

\begin{abstract} 
  This article provides entropic inequalities for binomial-Poisson
  distributions, derived from the two point space. They appear as local
  inequalities of the M/M/$\infty$ queue. They describe in particular the
  exponential dissipation of $\Phi$-entropies along this process. This simple
  queueing process appears as a model of ``constant curvature'', and plays for
  the simple Poisson process the role played by the Ornstein-Uhlenbeck process
  for Brownian Motion. Some of the inequalities are recovered by semi-group
  interpolation. Additionally, we explore the behaviour of these entropic
  inequalities under a particular scaling, which sees the Ornstein-Uhlenbeck
  process as a fluid limit of M/M/$\infty$
  queues. 
  Proofs are elementary and rely essentially on the development of a
  ``$\Phi$-calculus''.
\end{abstract}



\section{Introduction}

We consider in this article the M/M/$\infty$ queueing process. This elementary
continuous time Markov process on $\dN$ plays for the simple Poisson process
the role played by the Ornstein-Uhlenbeck process for Brownian motion. In
particular, its law at time $t$ is explicitly given by a binomial-Poisson
Mehler like formula, and the associated semi-group commutes with the discrete
gradient operator, up to a time decreasing exponential factor. We derive
general entropic inequalities for binomial-Poisson measures from the two
points space, essentially by convexity. They hold in particular for the law at
fixed time of the process, as for Ornstein-Uhlenbeck. In particular, these
entropic inequalities contain as special cases Poincar\'e inequalities and
various modified logarithmic Sobolev inequalities, which appear for instance
in \cite{MR1636948}, \cite{MR1757600}, \cite{MR1944012} and
\cite{boudou-caputo-daipra-posta}.

It is known that the lack of a chain rule and of a good notion of curvature in
discrete space settings make difficult the derivation of entropic inequalities
for discrete space Markov processes. Poincar\'e inequalities are exceptional,
due to their Hilbertian nature. Their derivation does not need the diffusion
property. L\'evy processes and Poisson space are also exceptional, since their
i.i.d. underlying structure makes them ``flat'' in a way. This nature is
translated on the infinitesimal Markov generator as a commutation with
translations. The M/M/$\infty$ queue has non-homogeneous independent
increments, and is thus beyond this framework. The reader may find various
entropic inequalities for finite space Markov processes in
\cite{MR1410112,MR1490046} and \cite{MR1815775}, and for infinite countable
space Markov processes in \cite{MR1710983}, \cite{MR1757600},
\cite{MR1800540}, \cite{MR1953566}, \cite{MR2023358}, \cite{MR2085608},
\cite{caputo-posta,boudou-caputo-daipra-posta}, \cite{MR2081075},
\cite{MR2184099,MR2147317,MR1944012}, \cite{joulin} for instance.

Birth and death processes are the discrete space analogue of diffusion
processes. However, they are not diffusions, and specific diffusion tools like
Bakry-\'Emery $\GD$ calculus are of difficult usage for such processes. It
follows from our study that the M/M/$\infty$ queueing process can serve as a
model of ``constant curvature'' on $\dN$. It is known that convexity may serve
as an alternative to the diffusion property, as presented for instance in
\cite{MR1800540} and \cite{MR2081075}. In this article, we circumvent the lack
of chain rule by elementary convexity bounds for germs of discrete Dirichlet
forms. This work can be seen as a continuation of \cite{MR2081075}, and was
initially motivated by the time inhomogeneous M/M/$\infty$ queue which appears
in the biological problem studied in \cite{chafai-concordet}. The notion of
queueing process is widely used in applied probability. The reader may find a
modern introduction to queueing processes in the book \cite{MR1996883} by
Robert, in connection with random networks, general Markov processes,
martingales, and point processes. This large family of Markov processes
contains, as particular cases, the simple Poisson process, the continuous time
simple random walk on $\dN$, and more generally all continuous time birth and
death processes on $\dN$.

The approach and results of this article may be extended by various ways. The
first step is to consider birth and death processes on $\dN$ or $\dZ$, and
then on $\dN^d$ or $\dZ^d$ with interactions. Some versions of such models
where already considered in the statistical mechanics literature, at least for
Poincar\'e and modified logarithmic Sobolev inequalities, see for instance
\cite{caputo-posta,boudou-caputo-daipra-posta} \cite{MR2184099,MR2147317},
\cite{MR1944012}, and references therein. These extensions concern continuous
time processes on a discrete space $E$ with generator of the form
$$
\bL(f)(x) = \int_E (f(y)-f(x))\,d\ga_x(y),
$$
where $\ga_x$ is the ``jump measure'' at point $x$, which is a finite Borel
measure on $E$. Another possibility is to consider Volterra processes driven
by a simple Poisson process, possibly together with Clark-Ocone formulas as in
\cite{MR1815775} and \cite{MR1800540} for instance. We hope that some of these
extensions will make the matter of forthcoming articles. We have in mind the
construction of a functional bridge between discrete space Markov processes
and ``curved'' diffusion processes, which complements, by mean of quantitative
functional inequalities, the approximation in law. The recent articles
\cite{MR2094147}, \cite{bobkov-tetali},
\cite{caputo-posta,boudou-caputo-daipra-posta},
\cite{MR2184099,MR2147317,MR1944012}, \cite{johnson-goldschmidt}, and
\cite{joulin} may help for such a program.

The entropic inequalities that we consider in this article can be called
``$\Phi$-Sobolev inequalities'' since they involve a $\Phi$-entropy and a
$\Phi$-Dirichlet form. They contain in particular Poincar\'e inequalities and
``$\rL^1$''-logarithmic-Sobolev inequalities. As presented in
\cite{MR2081075}, they hold, under convexity assumptions on $\Phi$, for
log-concave measures on $\dR^d$, for diffusions on manifolds with positive
bounded below curvature, for many Wiener measures, for Poisson space, and for
many L\'evy processes. Their genericity on $\Phi$ is particularly interesting
in discrete space settings for which no chain rule is available. The aim of
this article is to extend these entropic inequalities to discrete space
processes, beyond the i.i.d. increments case, in particular, beyond L\'evy
processes and Poisson space.

This work goes beyond many results of \cite{MR1636948}, in terms of entropies,
Dirichlet forms, and measures. We show how the entropic inequalities are
scaled when the discrete space curved process (M/M/$\infty$) approximates a
curved diffusion process (Ornstein-Uhlenbeck). This work can thus be seen as a
precise and instructive case study. Many aspects are still valid for more
general birth and death processes, and we believe that the entropic
inequalities that we consider here hold for interacting birth and death
processes. However, a lot of work is still needed to achieve this objective.
In particular, and to our knowledge, a good notion of curvature is still
lacking for interacting birth and death processes. Viewed as a unidimensional
(e.g. single site) particle system, the M/M/$\infty$ queue is not
conservative. In particular, it cannot be viewed as a particular
unidimensional case of the interacting birth and death processes considered in
\cite{caputo-posta,MR2184099,MR2147317}. However, it belongs to the models
considered in \cite{MR1944012}.

\textbf{Outline of the rest of the article.} In the introduction, the
definition of the M/M/$\infty$ queueing process is followed by the
presentation of links and analogies with the Ornstein-Uhlenbeck process, and
then by the introduction of the $\Phi$-entropy together with the $A-B-C$
transforms of $\Phi$. Section two is a two point space approach to
binomial-Poisson entropic inequalities. In Section three, we address the
exponential decay of $\Phi$-entropy functionals along the M/M/$\infty$ queue,
we give various proofs of entropic inequalities by using semi-group
interpolations, and we use a scaling limit to recover Gaussian inequalities
for the Ornstein-Uhlenbeck process. The fourth and last section is devoted to
key convexity properties related to the $A-B-C$ transforms.

\subsection{The M/M/$\infty$ queueing process}

The M/M/$\infty$ queue with input rate $\la\geq0$ and service rate $\mu\geq0$
is a particular space-inhomogeneous and time-homogeneous birth and death
process on $\dN$. Let $X_t$ be the number of customers in the queue -- i.e.
the length of the queue -- at time $t$. The name ``M/M/$\infty$'' comes from
Kendall's classification of simple queueing processes. It corresponds to an
infinite number of servers with random memoryless inter-arrivals (first M) and
service times (second M), see for instance \cite[p. xiii]{MR1996883}. Since
the number of servers is infinite, each client gets immediately his own
dedicated server, and the length of the queue is exactly the number of busy
servers. The infinitesimal Markov generator $\GI$ of the M/M/$\infty$ queue
$(X_t)_{t\geq0}$ is given for any $f:\dN\to\dR$ and any $n\in\dN$ by
\begin{equation}\label{eq:mmi-gi}
  \GI(f)(n) = n\mu\rD^*(f)(n)+\la\rD(f)(n),
\end{equation}
where the discrete gradients $\rD$ and $\rD^*$ are defined respectively for
any $f:\dZ\to\dR$ and any $n\in\dN$ by
\begin{equation}\label{eq:def-d-dstart}
  \rD(f)(n):=f(n+1)-f(n)
  \text{\quad and\quad}\rD^*(f)(n):=f(n-1)-f(n).
\end{equation}
The operators $\rD$ and $\rD^*$ commute with translations, but $\GI$ does not.
Notice that $f(-1)$ does not need to be defined in \eqref{eq:mmi-gi} since it
is multiplied by $0$. The stared notation for $\rD^*$ comes from the fact that
$\rD^*$ is the adjoint of $\rD$ with respect to the counting measure on $\dZ$.
The identity $\rD\rD^*=\rD^*\rD=-(\rD+\rD^*)$ leads to a polarised version of
the infinitesimal generator \eqref{eq:mmi-gi},
$$
\GI(f)(n)=-\la(\rD\rD^*)(f)(n)+(n\mu-\la)\rD^*(f)(n),
$$
for any $n\in\dN$ and $f:\dN\to\dR$. The finite difference operator $\rD\rD^*$
is the discrete Laplacian, given by $(\rD\rD^*)(f)(n) = 2f(n)-f(n-1)-f(n+1)$
for any $f:\dZ\to\dR$ and any $n\in\dZ$. 
Consider the process conditional to the event $\{X_s=n\}$. Let
$T:=\min\{t>s:X_t\neq X_s\}-s$ be the waiting time before next jump. Then $T$
follows an exponential law $\cE(\la+n\mu)$ of mean $1/(\la+n\mu)$. The
transition matrix $\bJ$ of the embedded discrete time jump Markov chain is
given for any $m,n\in\dN$ by
$$
\bJ(n,m):=
\frac{1}{\la+n\mu}
\begin{cases}
  n\mu & \text{if $m=n-1$} \\
  \la & \text{if $m=n+1$} \\
  0 & \text{otherwise}
\end{cases},
$$
where we assumed for simplicity that $\la+\mu>0$. The embedded chain is
recurrent irreducible as soon as $\la>0$ and $\mu>0$. The jump intensity
function $n\mapsto\la+n\mu$ is not bounded when $\mu>0$, however, the process
is not explosive by virtue of Reuter criterion for irreducible birth and
death processes, cf. \cite[Theorem 4.5]{MR1689633}.

Defining a stochastic process corresponds to specify a law on paths space.
Following \cite[Chap. 6]{MR1996883}, the stochastic process $(X_t)_{t\geq0}$
with $X_0=n$ can be constructed as follows~:
$$
X_t = n + \cN_\la(]0,t]) %
- \sum_{i=1}^\infty \int_{]0,t]}\!1_{\{X_{s^-}\geq i\}}\,\cN_\mu^i(ds),
$$
where $\cN_\la$ is a Poisson random measure on $\dR_+$ of intensity $\la$ and
where $(\cN_\mu^i)_{i\in\dN}$ is an i.i.d. collection of Poisson random
measures on $\dR_+$ of intensity $\mu$, independent of $\cN_\la$. In other
words, the process $(X_t)_{t\geq0}$ solves the Stochastic Differential
Equation
$$
dX_t = \cN_\la(dt) - \sum_{i=1}^{X_{t^-}} \cN^i_\mu(dt).
$$
Let us consider the filtration $(\cF_t)_{t\geq0}$ defined for any $t\in\dR_+$
by
$$
\cF_t:=\si\BRA{\cN_\la(]0,s]);s\in[0,t]} %
\vee\si\BRA{\cN_\mu^i(]0,s]);(s,i)\in[0,t]\times\dN}.
$$
The process $(X_t-X_0-\la t+\mu\int_0^t\!X_s\,ds)_{t\geq0}$ is a square
integrable martingale with increasing process given by $\la
t+\mu\int_0^t\!X_s\,ds$. More generally, the process $(X_t)_{t\geq0}$ is a
solution of the martingale problem associated to the Markov generator $\bL$
defined by \eqref{eq:mmi-gi}. Namely, for any $f:\dN\to\dR$, the process
$$
\PAR{f(X_t)-f(X_0)-\int_0^t\!\bL(f)(X_s)\,ds}_{t\geq0}
$$
is a local martingale. When $f(n)=n$ for any $n\in\dN$, we get
$\bL(f)(n)=\la-\mu n$. The Markov semi-group $(\SGf{t})_{t\geq0}$ of
$(X_t)_{t\geq0}$ is defined for any bounded $f:\dN\to\dR$ by
$$
\SG{t}{f}(n):=\moy{}{f(X_t)\,\vert X_0=n},
$$
in such a way that $\SG{t}{\rI_A}(n)=\dP(X_t\in A\,\vert\,X_0=n)$ for any
$A\subset\dN$. We have $\SG{t}{\cdot}(n)=\cL(X_t\,\vert X_0=n)$ for any
$n\in\dN$. In particular, $\SGf{t}\circ\SGf{s}=\SGf{t+s}$, and
$\SGf{0}=\mathrm{Id}$, and $\GI f:=\pd_{t=0}\SG{t}{f}$. The coefficient $\rho$
of the M/M/$\infty$ queue is defined by
$$
\rho:=\frac{\la}{\mu}.
$$
In the sequel, we denote by $\moy{Q}{f}$ or by $\moyf{Q}{f}$ the mean of
function $f$ with respect to the probability measure $Q$, and by $\rL^p(Q)$
the Lebesgue space of measurable real valued functions $f$ such that
$\ABS{f}^p$ is $Q$-integrable. For a Borel measure on $\dN$, we also denote
$Q(n):=Q(\{n\})$ for any $n\in\dN$.

\subsection{The Ornstein-Uhlenbeck process as a fluid limit of M/M/$\infty$ queues}
\label{ss:ou-mmi-fluid}

The Ornstein-Uhlenbeck process can be recovered from the M/M/$\infty$ queue as
a fluid limit, by using a Kelly scaling. See for instance
\cite{MR840104,MR1111523} and the books \cite{MR1996883}, \cite{MR838085}, and
\cite{MR1707314}. Namely, for any $N\in\dN$, let $(X^N_t)_{t\geq0}$ be the
M/M/$\infty$ queue with input rate $N\la$ and service rate $\mu>0$. Define the
process $(Y_t^N)_{t\geq0}$ by
$$
Y_t^N:=\frac{1}{N}X_t^N.
$$
For any $x\in\dR_+$, let $m:\dR_+\to\dR$ be defined by
$m(t):=\rho+(x-\rho)p(t)$ for any $t\in\dR_+$, where $p(t):=e^{-\mu t}$.
Consider a sequence of initial states $(X_0^N)_{N\in\dN}$ such that
$$
\lim_{N\to\infty} Y^N_0=\lim_{N\to\infty}\frac{1}{N}X_0^N=x.
$$
Then, for any $t\in\dR_+$, the sequence of random variables $\PAR{\sup_{0\leq
    s \leq t} \ABS{Y_s^N-m(s)}}_{N\in\dN}$ converges in $\rL^1$ towards $0$
when $N\to\infty$, see for instance \cite[Section 6.5]{MR1996883}. In
particular, for any $\veps>0$,
$$
\lim_{N\to\infty}\dP\PAR{\sup_{0\leq s\leq t}\ABS{Y_t^N-m(s)}>\veps} = 0.
$$
Moreover, this Law of Large Numbers is complemented by a Central Limit
Theorem, see for instance \cite{MR0222973} and \cite{MR1996883}. Namely,
define the process $(Z^N_t)_{t\geq0}$ by
$$
Z^N_t:=\sqrt{N}\PAR{Y^N_t-m(t)} = \frac{X^N_t-Nm(t)}{\sqrt{N}}.
$$
Notice that $m(0)=x$. Let $y\in\dR$ and assume that the initial states
$(X_0^N)_{N\in\dN}$ satisfy additionally that
$$
\lim_{N\to\infty} Z^N_0=\lim_{N\to\infty} \sqrt{N}\PAR{Y_0^N-x} = y.
$$
A basic example is given by $X_0^N=\PENT{Nx+\sqrt{N}y}$ where $\PENT{\cdot}$
denotes the integer part. Then, the sequence of processes $(Z^N_t)_{t\geq0}$
converges in distribution, when $N\to\infty$, towards a process denoted
$(Z^\infty_t)_{t\geq0}$, with non-homogeneous independent increments, given
by
$$
Z^\infty_t := y p(t) + \int_0^t p(t-s)\sqrt{\la+x\mu}\,dB_s,
$$
where $(B_s)_{s\geq0}$ is a standard Brownian Motion on the real line. In
particular, when $x=\rho$, then $m(s)=\rho$ for any $s\in\dR_+$, and
$(Z^\infty_t)_{t\geq0}$ is in that case an Ornstein-Uhlenbeck process,
solution of the Stochastic Differential Equation $Z^\infty_0=y$ and
$dZ^\infty_t=\sqrt{2\la}\,dB_t-\mu Z^\infty_tdt$, where $(B_t)_{t\geq0}$ is a
standard Brownian motion on the real line. Additionally, for any $t\in\dR_+$,
$$
\cL(Z^\infty_t\,\vert\,Z^\infty_0=y) =
\de_{yp(t)}*\cN\PAR{0,\PAR{1-p(t)^2}\rho},
$$
where $\cN(a,b)$ denotes the standard Gaussian law on $\dR$ of mean $a$ and
variance $b$. This Mehler formula is the continuous space analogue of
\eqref{eq:mehler}. The Markov infinitesimal generator of
$(Z^\infty_t)_{t\geq0}$ is the linear differential operator which maps
function $y\mapsto f(y)$ to function
$$
y\mapsto\la f''(y)-\mu yf'(y).
$$
The symmetric invariant measure of $(Z^\infty_t)_{t\geq0}$ is the Gaussian law
$\cN(0,\rho)$. The $\mu$ parameter appears clearly here as a curvature,
whereas the $\la$ parameter appears as a diffusive coefficient.

%
%

\subsection{The M/M/$\infty$ queue as a discrete Ornstein-Uhlenbeck process}

Let $(X_t)_{t\geq0}$ be an M/M/$\infty$ with rates $\la$ and $\mu$. When $\mu$
vanishes, $(X_t)_{t\geq0}$ reduces to a simple Poisson process of intensity
$\la$, and admits the counting measure on $\dN$ as a symmetric measure. A
contrario, when $\la$ vanishes, $(X_t)_{t\geq0}$ is a pure death process, and
admits $\de_0$ as an invariant probability measure.

The M/M/$\infty$ queue plays for the simple Poisson process the role played by
the Ornstein-Uhlenbeck process for standard Brownian Motion. The law of the
M/M/$\infty$ queue $(X_t)_{t\geq0}$ is explicitly given for any $n\in\dN$ by
the following Mehler like formula
\begin{equation}\label{eq:mehler}
  \cL(X_t\vert X_0=n)
  = \cB\PAR{n,p(t)}*\cP\PAR{\rho q(t)},
\end{equation}
where
$$
\quad p(t):=e^{-\mu t} \text{\quad and \quad} q(t):=1-p(t).
$$
When $\mu=0$, we set $\rho q(t)=\la$, since $\la=\lim_{\mu\to0^+}\rho q(t)$.
Here and in the sequel, $\cB(n,p)$ stands for the binomial distribution
$\cB(n,p):=(p\de_1+q\de_0)^{*n}$ of size $n\in\dN$ and parameter $p\in[0,1]$,
with the convention $\cB(n,0):=\de_0$ and $\cB(n,1):=\de_n$. The notation
$\cP(\si)$ stands for the Poisson measure on $\dN$ of intensity $\si>0$,
defined by $\cP(\si):=e^{-\si}\sum_{k=0}^{\infty}\frac{1}{k!}\si^k\de_k$. When
$\mu>0$, the process $(X_t)_{t\geq0}$ is ergodic and admits $\cP(\rho)$ as a
reversible invariant measure. In other words, for any $n\in\dN$ and
$s\in\dR_+$,
$$
\lim_{n\to\infty}\cL(X_t\vert X_s=n) = \cP(\rho).
$$
Moreover, $\moy{\cP(\rho)}{\SGf{t}{f}}=\moy{\cP(\rho)}{f}$ for any
$f\in\rL^1(\cP(\rho))$ and any $t\in\dR_+$, or equivalently
$\moy{\cP(\rho)}{\GI f}=0$ for any $f\in\rL^1(\cP(\rho))$. As for the
Ornstein-Uhlenbeck process, this convergence is not uniform in $n$ since for
any $\al>0$,
$$
\lim_{n\to\infty}\cL(X_{\mu^{-1}\log(n/\al)}\,\vert\,X_s=n)=\cP(\al+\rho).
$$
The mean and variance of $\cL(X_t\vert X_s=n)$ with $t\geq s\geq0$ are given
respectively by
$$
np(t-s) + \rho q(t-s) \text{\quad and \quad} (np(t-s)+\rho)q(t-s).
$$
The semi-group $(\SGf{t})_{t\geq0}$ of the M/M/$\infty$
queue shares the nice ``constant curvature'' property with the
Ornstein-Uhlenbeck semi-group. Namely, for any $t\in\dR_+$, any $n\in\dN$, and
any bounded $f:\dN\to\dR$,
\begin{equation}\label{eq:commut}
  \rD\SGf{t}f=e^{-\mu t}\SGf{t}\rD f.
\end{equation}
The infinitesimal version writes $[\bL,\rD]:=\bL\rD-\rD\bL=\mu\rD$. The
commutation \eqref{eq:commut} can be deduced simply from \eqref{eq:mehler}.
Namely, if $X_1,\ldots,X_{n+1},Y$ are independent random variables with
$X_i\sim\cB(1,p(t))$ and $Y\sim\cP(\rho q(t))$,
\begin{align*}
  \SG{t}{f}(n+1)
  & = \moy{}{f(X_1+\cdots+X_{n+1}+Y)} \\
  & = \moy{}{\moy{}{f(X_1+\cdots+X_{n+1}+Y)\vert X_{n+1}}} \\
  & = p(t)\SG{t}{f(1+\cdot)}(n)+q(t)\SG{t}{f}(n) \\
  & = p(t)\SG{t}{\rD f}(n)+\SG{t}{f}(n).
\end{align*}
This fact and the properties of the $A-B-C$ transforms introduced in the
sequel give rise to various entropic inequalities, by using the semi-group
$(\SGf{t})_{t\geq0}$ as an interpolation flow.

We give now various binomial-Poisson ``integration by parts'' formulas. Let
$H_{n,p}(m):=\binom{n}{m}p^mq^{n-m}$ for any $p\in[0,1]$ and any integers $n$
and $m$ with $0\leq m\leq n$. We have then $mH_{n,p}(m)=npH_{n-1,p}(m-1)$ as
soon as $0<m\leq n$. As a consequence, for any function $f:\dN\to\dR$ , any
$n\in\dN^*$ and any $p\in[0,1]$
\begin{equation}\label{eq:ipp-bin}
  \moy{\cB(n,p)}{hf}
  =np\moy{\cB(n-1,p)}{f(1+\cdot)},
\end{equation}
where $h:\dN\to\dR$ is defined by $h(k)=k$ for any $k\in\dN$. Similarly,
$(n-m)H_{n,p}(m)=nqH_{n-1,p}(m)$ as soon as $0\leq m<n$, which gives
for any function $f:\dN\to\dR$ , any $n\in\dN^*$ and any $p\in[0,1]$
\begin{equation}\label{eq:ipp-bin-bw}
  \moy{\cB(n,p)}{(n-h)f}
  =nq\moy{\cB(n-1,p)}{f}.
\end{equation}
For $\rho>0$, the binomial approximation of Poisson measure which lets $np$
tend to $\rho$ when $n\to\infty$ gives from \eqref{eq:ipp-bin}
\begin{equation}\label{eq:ipp-poi}
  \moy{\cP(\rho)}{hf}
  =\rho\moy{\cP(\rho)}{f(1+\cdot)}.
\end{equation}
Some algebra provides, by conditioning, a mixed binomial-Poisson version
\begin{equation}\label{eq:ipp-bin-poi}
  \moy{\cB(n,p)*\cP(\rho)}{hf}
  =np\moy{\cB(n-1,p)*\cP(\rho)}{f(1+\cdot)}
  +\rho\moy{\cB(n,p)*\cP(\rho)}{f(1+\cdot)}.
\end{equation}
In particular, the Mehler like formula \eqref{eq:mehler} gives for any
$n\in\dN^*$ and $t\in\dR_+$,
\begin{equation}\label{eq:ipp-bin-poi-sg}
  \mu\SG{t}{hf}(n) %
  = \mu n p(t)\SG{t}{f(1+\cdot)}(n-1) + \la q(t)\SG{t}{f(1+\cdot)}(n),
\end{equation}
where $h:\dN\to\dN$ is defined by $h(n):=n$ for any $n\in\dN$. The
binomial-Poisson nature of the M/M/$\infty$ queue is related to the fact that
the coefficients of its infinitesimal generator \eqref{eq:mmi-gi} are affine
functions of $n$. The reader may find an analysis of linear growth birth and
death processes in \cite{MR0098435} and \cite{MR1056442} and references
therein.

\subsection{Convex functionals}

For any convex domain $\cD$ of $\dR^n$, let us denote by $\cC_\cD$ the convex
set of smooth convex functions from $\cD$ to $\dR$. Let $\cI\subset\dR$ be an
open interval of $\dR$ and $\Phi\in\cC_\cI$. We denote by $\rL^{1,\Phi}(Q)$
the convex subset of functions $f\in\rL^1(Q)$ taking their values in $\cI$ and
such that $\Phi(f)\in\rL^1(Q)$. We define the $\Phi$-entropy
$\ENT{Q}{\Phi}{f}$ of $f\in\rL^{1,\Phi}(Q)$ by
$$
\ENT{Q}{\Phi}{f}:=\moy{Q}{\Phi(f)}-\Phi(\moyf{Q}{f}).
$$
Is is also known as ``Jensen divergence'' since Jensen inequality gives
$\ENT{Q}{\Phi}{f}\geq0$. Moreover, when $\Phi$ is strictly convex,
$\ENT{Q}{\Phi}{f}=0$ if and only if $\Phi(f)$ is $Q$-a.s. constant. One can
distinguish for function $\Phi$ the following three usual special cases.
\begin{enumerate}
\item[\textbf{(P1)}]\label{it:phi-case-xlogx} $\Phi(u)=u\log(u)$ on
  $\cI=\dR_+^*$, and $\ENT{Q}{\Phi}{f}=\ent{Q}{f}:=\moy{Q}{f\log
    (f/\moyf{Q}{f})}$;
\item[\textbf{(P2)}]\label{it:phi-case-xx} $\Phi(u)=u^2$ on $\cI=\dR$, and
  $\ENT{Q}{\Phi}{f}=\var{Q}{f}:=\moy{Q}{(f-\moyf{Q}{f})^2}$;
\item[\textbf{(P3)}]\label{it:phi-case-xalpha} $\Phi(u)=u^\al$ on $\cI=\dR_+^*$
  with $\al\in(1,2)$.
\end{enumerate}
The $\ENTF{Q}{\Phi}$ functional is linear in $\Phi$ and vanishes when $\Phi$
is affine. Let us define from the interval $\cI\subset\dR$ the convex subsets
$\cT_\cI\subset\cT'_\cI\subset\dR^2$ by
$$
\cT_\cI:=\BRA{(u,v)\in\dR^2;\,(u,u+v)\in\cI\times\cI}
\text{\quad and\quad}
\cT'_\cI:=\BRA{(u,v)\in\dR^2;\,u\in\cI,(v+\cI)\cap\cI\neq\emptyset}.
$$
The $A-B-C$ transforms of $\Phi$ are the functions
$A^\Phi,B^\Phi:\cT_\cI\to\dR$ and $C^\Phi:\cT'_\cI\to\dR$ defined by
\begin{eqnarray*}
  A^\Phi(u,v) &:= & \Phi(u+v)-\Phi(u)-\Phi'(u)v; \\
  B^\Phi(u,v) &:= & (\Phi'(u+v)-\Phi'(u))v; \\
  C^\Phi(u,v) &:= & \Phi''(u)v^2.
\end{eqnarray*}
These three transforms are linear in $\Phi$, and their kernel contains any
affine function. Various additional properties of these three transforms are
collected in Section \ref{se:cvx}. In particular, the convexity of $\Phi$ on
$\cI$ is equivalent to the non negativity of its $A-B-C$ transforms on
$\cT_\cI$. In particular, the following statements hold.
\begin{itemize}
\item $A^\Phi$, $B^\Phi$, $C^\Phi$, $\ENTF{Q}{\Phi}$ are non negative and
  convex for \textbf{(P1-P2-P3)};
\item $2A^\Phi=B^\Phi=C^\Phi$ for \textbf{(P2)}, $A^\Phi\leq C^\Phi$ for
  \textbf{(P1)}, and $A^\Phi\leq B^\Phi$ for \textbf{(P1-P2-P3)}.
\end{itemize}
On the two point space $\{0,1\}$, the $\Phi$-entropy gives rise naturally to
the $A$-transform of $\Phi$. Namely, for any $f:\{0,1\}\to\cI$
with $(a,b):=(f(0),f(1))$ and $(u,v):=(a,b-a)$,
\begin{equation}\label{eq:ent-twop}
  \ENT{\cB(1,p)}{\Phi}{f}
  =q\Phi(a)+p\Phi(b)-\Phi(qa+pb) 
  =pA^\Phi(u,v)-A^\Phi(u,pv).
\end{equation}
The $A-B-C$ transforms are the germs of discrete Dirichlet forms via the
identities
\begin{align*}
  A^\Phi(f,\rD f)
  &=\rD\Phi(f)-\Phi'(f)\rD f; \\
  B^\Phi(f,\rD f)
  &=\rD(\Phi'(f))\rD f; \\
  C^\Phi(f,\rD f)
  &=\Phi''(f)\ABS{\rD f}^2.
\end{align*}
The reader may find explicit examples in table \ref{tab:ABCxpl}. We used above
the following identity, valid for any functions $\vphi:\dR\to\dR$ and
$f:\dZ\to\dR$,
$$
\rD(\vphi(f)) = \vphi(f(1+\cdot)) - \vphi(f) = \vphi(f+\rD f) - \vphi(f),
$$
where $f(1+\cdot)$ stands for $\dZ\ni n\mapsto f(1+n)$. In particular,
$f(1+\cdot)=f+\rD f$. The usage of the $A-B-C$ transforms allows, as presented
in the sequel, to derive several entropic inequalities in the same time,
including Poincar\'e inequalities and various modified logarithmic Sobolev
inequalities. They reduce most of the proofs to convexity, and they provide
various comparisons for discrete Dirichlet forms.

For any open interval $\cI\subset\dR$ and any probability space $(E,\cE,Q)$,
we denote in the sequel by $\cK(E,\cI)$ the convex set of measurable functions
$E\to\cI$ with a relatively compact image in $\cI$. These functions are
bounded. Notice that $\cK(E,\cI)$ is a convex subset of $\rL^{1,\Phi}(Q)$. The
introduction of $\cK(E,\cI)$ permits to avoid integrability obstructions at
the boundary of $\cI$ when dealing with the derivatives of $\Phi$. Any element
of $\rL^{1,\Phi}(Q)$ is a pointwise limit of elements of $\cK(E,\cI)$.

\begin{table}[htbp]
  \begin{center}
    \begin{tabular}{l|l|l|l} 
      \textbf{Function} $\BS{\Phi}$
      & $\BS{\cI}$
      & $\BS{A^\Phi}$
      & $\BS{A^\Phi(f,\rD f)}$ \\ [1em] 
      \textbf{(P1)} $u\log(u)$
      & $\dR_+^*$
      & $\ds{(u+v)(\log(u+v)-\log(u))-v}$
      & $(f+\rD f)\rD(\log f)-\rD f$\\
      \textbf{(P2)} $u^2$
      & $\dR$
      & $v^2$
      & $\ABS{\rD f}^2$ \\
      \textbf{(P3)} $u^\al$
      & $\dR_+^*$
      & $\ds{(u+v)^\al-u^\al-\al u^{\al-1}v}$
      & $\rD(f^\al)-\al f^{\al-1}\rD f$ \\ [1em]
      &
      & $\BS{B^\Phi}$
      & $\BS{B^\Phi(f,\rD f)}$ \\ [0.5em] 
      \textbf{(P1)} $u\log(u)$
      & $\dR_+^*$
      & $\ds{v(\log(u+v)-\log(u))}$
      & $\rD(f)\rD(\log f)$ \\
      \textbf{(P2)} $u^2$
      & $\dR$
      & $2v^2$
      & $2\ABS{\rD f}^2$ \\
      \textbf{(P3)} $u^\al$
      & $\dR_+^*$
      & $\ds{\al v((u+v)^{\al-1}-u^{\al-1})}$
      & $\al \rD(f)\rD(f^{\al-1})$  \\ [1em]
      &
      & $\BS{C^\Phi}$
      & $\BS{C^\Phi(f,\rD f)}$ \\ [0.5em] 
      \textbf{(P1)} $u\log(u)$
      & $\dR_+^*$
      & $\ds{v^2 u^{-1}}$
      & $\ABS{\rD f}^2 f^{-1}$ \\
      \textbf{(P2)} $u^2$
      & $\dR$
      & $2v^2$
      & $2\ABS{\rD f}^2$ \\
      \textbf{(P3)} $u^\al$
      & $\dR_+^*$
      & $\ds{\al(\al-1)v^2u^{\al-2}}$
      & $\al(\al-1)\ABS{\rD f}^2 f^{\al-2}$ 
    \end{tabular}
    \caption{Examples of $A-B-C$ transforms. For \textbf{(P3)},
      $\al\in(1,2)$.}
  \label{tab:ABCxpl}
\end{center}
\end{table}

\section{From two point space to binomial-Poisson inequalities}

Let $p\in[0,1]$ and let $\cB(1,p)$ be the Bernoulli measure $q\de_0+p\de_1$ on
$\{0,1\}$, where $q:=1-p$. We identify the two point space $\{0,1\}$ with the
``circle'' $\dZ/2\dZ$, for which $1+1=0$. In particular, the the ``+'' sign in
the definition \eqref{eq:def-d-dstart} of $\rD$ is taken modulo $2$. Then, for
any $f:\{0,1\}\to\cI$, the following identity holds.
$$
pq\moy{\cB(1,p)}{B^\Phi(f,\rD f)}- \ENT{\cB(1,p)}{\Phi}{f} %
=A^\Phi(\si_p(a,b-a))+A^\Phi(\si_q(b,a-b)),
$$
where $(a,b):=(f(0),f(1))$ and where $\si_p$ is as in \eqref{eq:def-sigp}.
Now, Lemma \ref{le:cvx-1} gives that $A^\Phi$ is non-negative as soon as
$\Phi$ is convex. Consequently, when $\Phi$ is convex, we obtain the following
entropic inequality for $\cB(1,p)$.
\begin{equation}\label{eq:bern-mlsi-B}
  \ENT{\cB(1,p)}{\Phi}{f} %
  \leq pq\moy{\cB(1,p)}{B^\Phi(f,\rD f)}.
\end{equation}
Unfortunately, the inequality \eqref{eq:bern-mlsi-B} is not optimal for
\textbf{(P2)} since in that case,
$$
\ENT{\cB(1,p)}{\Phi}{f}=\var{\cB(1,p)}{f}=pq(f(1)-f(0))^2 %
\text{\quad whereas\quad} %
\moy{\cB(1,p)}{B^\Phi(f,\rD f)}=2(f(1)-f(0))^2.
$$
This is due to the fact that $B^\Phi(f,\rD f)=2\ABS{\rD f}^2$ for
\textbf{(P2)}. We derive in the sequel the $A$ transform version, which is
stronger and optimal for \textbf{(P2)} since $A^\Phi(f,\rD f)=\ABS{\rD f}^2$
in that case. All the inequalities obtained in this section involve the $A$
transform in their right hand side. They hold for example in the cases
\textbf{(P1)}, \textbf{(P2)}, \textbf{(P3)}. The $A$ transform can be bounded
by the $B$ or the $C$ transforms, by using the elementary bounds given by
Lemma \ref{le:cvx-3/2}. We start with an entropic inequality for the Bernoulli
law $\cB(1,p)$. By convolution, we derive from this two point space
inequality a new entropic inequality for the binomial law
$\cB(n,p)=\cB(1,p)^{*n}$. An inequality for the Poisson law $\cP(\rho)$ is
then obtained by binomial approximation. The binomial-Poisson case is derived
by tensorisation. The following calculus Lemma is a $\Phi$ version of
\cite[Lemma 2]{MR1636948} by Bobkov and Ledoux.

\begin{lem}[two point Lemma]\label{le:cond-u}
  Let $\Phi\in\cC_\cI$ such that $\Phi''\in\cC_\cI$. Let $U:[0,1]\to\dR$ be
  defined by
  \begin{equation*}
    U(p):=\ENT{\cB(1,p)}{\Phi}{f}-pq\moy{\cB(1,p)}{g},
  \end{equation*}
  where $f,g:\{0,1\}\to\cI$. Then, $U\leq0$ on $[0,1]$ if and
  only if
  \begin{equation}\label{eq:cond-u}
    U'(0)\leq 0\leq U'(1).
  \end{equation}
\end{lem}

\begin{proof}
  We denote $(a,b):=(f(0),f(1))$ and $(\al,\be):=(g(0),g(1))$. We get then
  \begin{equation}\label{eq:defU}
    U(p)=q\Phi(a)+p\Phi(b)-\Phi(qa+pb)-pq(q\al+p\be).
  \end{equation}
  The last term is a polynomial in $p$ of degree three. Taking the fourth
  derivative in $p$ gives
  \begin{equation*}
    U''''(p)=-(b-a)^4\Phi''''(qa+pb).
  \end{equation*}
  Since $\Phi''$ is convex, we have $U''''\leq0$ on $(0,1)$ and thus $U''$ is
  concave. Consequently, there exists $0\leq p_0 \leq p_1 \leq 1$ such that
  $U''\leq0$ on $[0,p_0]\cup[p_1,1]$ and $U''\geq0$ on $[p_0,p_1]$. We have
  that $U$ is concave on $[0,p_0]$. But $U(0)=0$ and by assumption
  $U'(0)\leq0$, thus $U\leq0$ on $[0,p_0]$ by concavity. A consequence is that
  $U(p_0)\leq0$. Similarly by symmetry we have that $U\leq0$ on $[p_1,1]$ and
  $U(p_1)\leq0$. Now since $U$ is convex on $[p_0,p_1]$ and non-positive on
  the boundaries, it is non-positive on the whole interval $[p_0,p_1]$.
  Therefore \eqref{eq:cond-u} implies $U\leq0$ on $[0,1]$.
\end{proof}

One can show by similar arguments that if additionally $f(0)\geq f(1)$ and
$g(0)\geq g(1)$ (respectively $f(0)\leq f(1)$ and $g(0)\leq g(1)$), then the
condition \eqref{eq:cond-u} may be weakened into $U'(0)\leq0$ (respectively
$U'(1)\geq0$). Notice that in terms of $A$ transform,
\begin{equation}\label{eq:uprime-0-1}
  U'(0)=A^\Phi(a,b-a)-\al
  \text{\quad and \quad}
  U'(1)=A^\Phi(a,b-a)+\be.
\end{equation}

The following Lemma provides the $A$ transform version of
\eqref{eq:bern-mlsi-B}.

\begin{lem}[two point entropic inequalities]\label{le:twop}
  Let $\Phi\in\cC_\cI$ such that $\Phi''\in\cC_\cI$. Then, for any
  $f:\{0,1\}\to\cI$,
  \begin{equation}\label{eq:bern-mlsi-A}
    \ENT{\cB(1,p)}{\Phi}{f}
    \leq pq\,\moy{\cB(1,p)}{A^\Phi(f,\rD f)},
  \end{equation}
  where the ``+'' in \eqref{eq:def-d-dstart} of $\rD$ is taken modulo $2$.
  Moreover, the inequality becomes an equality for \textbf{(P2)}.
\end{lem}

\begin{proof}
  Let $U$ be as in \eqref{eq:defU} with $g=A^\Phi(f,\rD f)$. From
  \eqref{eq:uprime-0-1} we get
  \begin{equation*}
    U'(0)=A^\Phi(a,b-a)-A^\Phi(a,b-a)=0
    \text{\quad and \quad}
    U'(1)=A^\Phi(a,b-a)+A^\Phi(b,a-b),
  \end{equation*}
  where $(a,b):=(f(0),f(1))$. By Lemma \ref{le:cvx-1}, function $A^\Phi$ is
  non negative and thus $U'(1)\geq0$. Therefore \eqref{eq:bern-mlsi-A} follows
  by virtue of Lemma \ref{le:cond-u}. Notice that since $+1=-1$ in $\dZ/2\dZ$,
  we have $\rD=\rD^*$. In particular, for any $f:\{0,1\}\to\cI$,
  the function $B^\Phi(f,\rD f)$ is constant, and $A^\Phi(f,\rD
  f)=A^\Phi(f,\rD^* f)$.
\end{proof}

Notice that \eqref{eq:bern-mlsi-A} can be rewritten as \eqref{eq:pA-Ap}.
Entropic inequalities like \eqref{eq:bern-mlsi-A} belong to the so called
family of $\Phi$-Sobolev inequalities, which are known to be stable by
convolution, cf. \cite[Corollary 3.1 page 342]{MR2081075}. This observation
leads to Theorem \ref{th:gbern} below, by using the tensorisation property
\eqref{eq:tenso} of Theorem \ref{th:cvx-2}.

\begin{thm}[Bernoulli entropic inequalities]\label{th:gbern}
  Let $M:=\cB(1,p_1)*\cdots*\cB(1,p_n)$ and $C_M:=\max\{p_1q_1,\ldots,p_n
  q_n\}$ where $p_1,\ldots,p_n\in[0,1]$. Let $\Phi\in\cC_\cI$ such that
  $A^\Phi\in\cC_{\cT_\cI}$. Then, for any $f:\dN\to\cI$,
  \begin{equation}\label{eq:gbern-A}
    \ENT{M}{\Phi}{f}
    \leq C_M \moy{M}{(n-h)A^\Phi(f,\rD f)+h A^\Phi(f,\rD^* f)},
  \end{equation}
  where $h:\dN\to\dR$ is defined by $h(k)=k$ for any $k\in\dN$. In particular,
  \begin{equation}\label{eq:bin-A}
    \ENT{\cB(n,p)}{\Phi}{f}
    \leq pq \moy{\cB(n,p)}{(n-h)A^\Phi(f,\rD f)+h A^\Phi(f,\rD^* f)},
  \end{equation}
  for any $n\in\dN^*$, any $p\in[0,1]$, and any $f:\dN\to\cI$.
  Moreover, if $\tau$ is as in \eqref{eq:def-tau},
  \begin{equation}\label{eq:bin-Abis}
    \ENT{\cB(n,p)}{\Phi}{f} \leq %
    npq \moy{\cB(n-1,p)}{qA^\Phi(f,\rD f)+pA^\Phi(\tau(f,\rD f))}.
  \end{equation}
\end{thm}

The optimality of these inequalities in the case \textbf{(P2)} can be checked
for a linear function $f$.

\begin{proof}
  First of all, by virtue of Theorem \ref{th:cvx-2}, the convexity of $A^\Phi$
  on $\cT_\cI$ implies the convexity of $\Phi''$ on $\cI$. Let
  $(E_i,Q_i)=(\{0,1\},\cB(1,p_i))$ for any $i\in\{1,\ldots,n\}$. Let
  $f:\dN\to\cI$ and consider the symmetric function
  $g:E_1\times\cdots\times E_n\to\cI$ defined by
  $g(x_1,\ldots,x_n):=f(x_1+\cdots+x_n)$. The tensorisation formula
  \eqref{eq:tenso} together with the two point entropic inequality
  \eqref{eq:bern-mlsi-A} of Lemma \ref{le:twop} gives
  \begin{equation*}
    \ENT{Q_1\otimes\cdots\otimes Q_n}{\Phi}{g}
    \leq C_Q
    \moy{Q_1\otimes\cdots\otimes Q_n}{\sum_{i=1}^n A^\Phi(g,\rD^{(i)} g)},
  \end{equation*}
  where $\rD^{(i)}$ denotes the operator $\rD$ acting on the $i^\text{th}$
  coordinate with modulo $2$ as in Lemma \ref{le:twop}. At this step, we
  notice by denoting $s_n:=x_1+\cdots+x_n$ that for any $x\in\{0,1\}^n$,
  $$
  \sum_{i=1}^n A^\Phi(g,\rD^{(i)} g)(x) %
  =(n-s_n)A^\Phi(f,\rD f)(s_n)+s_nA^\Phi(f,\rD^* f)(s_n).
  $$
  Outside $\{0,\ldots,n\}$, the function $f$ takes values which come with a
  null coefficient in the right hand side. The desired result follows since
  $M$ is the law of $s_n$ under $Q_1\otimes\cdots\otimes Q_n$. Inequality
  \eqref{eq:gbern-A} reduces to \eqref{eq:bin-A} when $p_1=\cdots=p_n=p$. It
  remains to establish \eqref{eq:bin-Abis}. By virtue of \eqref{eq:ipp-bin}
  and \eqref{eq:ipp-bin-bw}, the right hand side of \eqref{eq:bin-A} is equal
  to
  $$
  n\moy{\cB(n-1,p)}{qA^\Phi(f,\rD f)+pA^\Phi(f,\rD^* f)(1+\cdot)}.
  $$
  Inequality \eqref{eq:bin-Abis} follows then from the simple identity
  \begin{equation}\label{eq:adtau}
   (f,\rD^* f)(1+\cdot)=\tau(f,\rD f).
  \end{equation}
  When $n=1$, then $M=\cB(1,p)$, and \eqref{eq:gbern-A} reduces to
  \eqref{eq:bern-mlsi-A}. Beware that $\rD$ in \eqref{eq:bern-mlsi-A} is taken
  modulo $2$.
\end{proof}

\begin{cor}[Poisson entropic inequality]\label{co:pois}
  Let $\Phi\in\cC_\cI$ be such that $A^\Phi\in\cC_{\cT_\cI}$. Let $\rho>0$ and
  $\cP(\rho)$ be the Poisson measure of mean $\rho$. Then, for any
  $\rho\in\dR_+$ and any $f\in\rL^{1,\Phi}(\cP(\rho))$,
  \begin{equation}\label{eq:pois-A}
    \ENT{\cP(\rho)}{\Phi}{f}
    \leq \rho\,\moy{\cP(\rho)}{A^\Phi(f,\rD f)}.
  \end{equation}
\end{cor}

\begin{proof}
  Notice that the right hand side of \eqref{eq:pois-A} takes its values in
  $[0,+\infty]$. By approximation, we can assume that $f\in\cK(\dN,\cI)$.
  Consider now \eqref{eq:bin-Abis}. Let $p$ depend on $n$ is such a way that
  $\lim_{n\to\infty}np_n=\rho$. Since $\lim_{n\to\infty}p_n=0$, we have
  $q_n\to1$. Moreover, $\cB(n,p_n)\to\cP(\rho)$ and
  $\cB(n-1,p_n)\to\cP(\rho)$.
\end{proof}

To the author's knowledge, inequality \eqref{eq:pois-A} appeared for the first
time in \cite{MR1800540} for \textbf{(P1)}, in \cite[p. 357]{MR1636948} for
\textbf{(P2)}, and in \cite{MR2081075} for the general case. See also
\cite{boudou-caputo-daipra-posta}. The $B$ and $C$ transforms versions of
\eqref{eq:pois-A}, which are weaker, appeared in particular in
\cite{MR1757600} and \cite{MR1636948}.

\begin{cor}[Binomial-Poisson entropic inequality]\label{co:binpois}
  Let $\Phi\in\cC_\cI$ be such that $A^\Phi\in\cC_{\cT_\cI}$. Let $M_n$ be the
  probability measure $M_n=\cB(n,p)*\cP(\rho)$ where $p\in[0,1]$,
  $\rho\in\dR_+$, and $n\in\dN$. Then, for any $f\in\rL^{1,\Phi}(\cP(\rho))$,
  \begin{equation}\label{eq:binpois-A}
    \ENT{M_n}{\Phi}{f} %
    \leq %
    \rho\moy{M_n}{A^\Phi(f,\rD f)} %
    + npq\moy{M_{n-1}}{qA^\Phi(f,\rD f)+pA^\Phi(\tau(f,\rD f))}.
  \end{equation}
\end{cor}

\begin{proof}
  By approximation, we can assume that $f\in\cK(\dN,\cI)$. If $n=0$, then
  \eqref{eq:binpois-A} reduces to \eqref{eq:pois-A}. Let us assume now that
  $n>0$. Let $(E_1,Q_1)=(\dN,\cB(n,p))$ and $(E_2,Q_2)=(\dN,\cP(\rho))$. Let
  $g:E_1\times E_2\to\cI$ be defined by $g(x_1,x_2)=f(x_1+x_2)$.
  Let us denote by $\rD^{(1)}$ and $\rD^{(2)}$ the $\rD$ operator which acts
  on $x_1$ and $x_2$ respectively. The inequalities \eqref{eq:tenso},
  \eqref{eq:bin-Abis}, \eqref{eq:pois-A} yield that $\ENT{Q_1\otimes
    Q_2}{\Phi}{g}$ is bounded above by
  $$
  npq\moy{Q_2}{\moy{Q_0}{qA^\Phi(g,\rD^{(1)} g)+pA^\Phi(\tau(g,\rD^{(1)} g)})}
  +\rho\moy{Q_1}{\moy{Q_2}{A^\Phi(g,\rD^{(2)} g)}},
  $$
  where $Q_0:=\cB(n-1,p)$. Since $\rD$ commutes with translations, we get
  for $i=1,2$,
  $$
  (g,\rD^{(i)}g)(x_1,x_2)=(f,\rD f)(x_1+x_2).
  $$
  Inequality \eqref{eq:binpois-A} follows since $M_n$, respectively $M_{n-1}$,
  is the law of $x_1+x_2$ under $Q_1\otimes Q_2$, respectively $Q_0\otimes
  Q_2$.
\end{proof}

The expectation with respect to $M_{n-1}$ in the right hand side of
\eqref{eq:binpois-A} may be rewritten as an expectation with respect to $M_n$
by using \eqref{eq:ipp-bin-poi}.

\section{Entropies along the M/M/$\infty$ queue}

We start with the decay of the $\Phi$-entropy functional along the queue.

\begin{thm}[$\Phi$-entropies dissipation]\label{th:MLSI}
  Let $\Phi\in\cC_\cI$. Let $(\SGf{t})_{t\geq0}$ be the M/M/$\infty$
  semi-group with input rate $\la>0$ and service rate $\mu>0$. Then for any
  $f\in\cK(\dN,\cI)$, the real function
  $t\in\dR_+\mapsto\ENT{\cP(\rho)}{\Phi}{\SGf{t}{f}}$ is non-increasing.
  Moreover, when $A^\Phi\in\cC_{\cT_\cI}$,
  \begin{equation}\label{eq:entropy-dec}
   \ENT{\cP(\rho)}{\Phi}{\SGf{t}{f}} %
   \leq e^{-c \mu t} \ENT{\cP(\rho)}{\Phi}{f},
  \end{equation}
  where $c$ is the best (i.e. biggest) constant in the inequality
  $$
  \forall f\in\rL^{1,\Phi}(\cP(\rho)),\quad %
  c\mu\ENT{\cP(\rho)}{\Phi}{f}\leq \la\moy{\cP(\rho)}{B^\Phi(f,\rD f)}.
  $$
  It holds with $c=1$ in general, and with $c=2$ for \textbf{(P2)}.
\end{thm}

\begin{proof}
  Let us denote $Q=\cP(\rho)$. Since
  $\ENT{Q}{\Phi}{\SGf{t}{f}}=\moy{Q}{\SGf{t}{\Phi(f)}}-\Phi(\moy{Q}{f})$, the
  invariance of $Q$ gives,
  $$
  \pd_t \ENT{Q}{\Phi}{\SGf{t}{f}}=\moy{Q}{\Phi'(\SGf{t}f)\GI\SGf{t} f}.
  $$
  Jensen inequality yields $\SG{t}{\Phi(f)}\geq\Phi(\SG{t}{f})$ as soon as
  $\Phi$ is convex. In particular $\GI\Phi(f)\geq\Phi'(f)\GI f$, which gives
  $\moy{Q}{\Phi'(f)\GI f}\leq0$ as soon as $\Phi(f)$ is $Q$-integrable. Hence,
  $t\mapsto\ENT{Q}{\Phi}{\SGf{t}{f}}$ is non-increasing, and we used only the
  convexity of $\Phi$, the Markovian nature of $(\SGf{t})_{t\geq0}$, and the
  invariance of $Q$. The Poisson integration by parts \eqref{eq:ipp-poi} --
  which is this time specific to our settings -- yields for any $g$
  \begin{equation}\label{eq:propBPoi}
    \moy{Q}{\Phi'(g)\GI g} %
    = -\la\moy{Q}{\rD(g)\rD(\Phi'(g))} %
    = -\la\moy{Q}{B^\Phi(g,\rD g)}.
  \end{equation}
  In particular, for $g=\SG{t}{f}$, we get,
  $$
  \pd_t\ENT{Q}{\Phi}{\SGf{t}{f}} =
  -\la\moy{Q}{B^\Phi(\SGf{t}{f},\rD\SGf{t}{f})}.
  $$
  Notice that since $\Phi$ is convex, we have $B^\Phi\geq0$ by virtue of Lemma
  \ref{le:cvx-1}. In the other hand, when $A^\Phi$ is convex, the Poisson
  entropic inequality \eqref{eq:pois-A} together with the bound $A^\Phi\leq
  B^\Phi$ of Lemma \ref{le:cvx-3/2} gives
  $$
  -\la\moy{Q}{B^\Phi(\SGf{t}f,\rD \SGf{t}f)} %
  \leq -\mu\ENT{Q}{\Phi}{\SGf{t}{f}}.
  $$
  Putting all together yields $\pd_t\ENT{Q}{\Phi}{\SGf{t}{f}}
  \leq-\mu\ENT{Q}{\Phi}{\SGf{t}{f}}$, which gives \eqref{eq:entropy-dec}.
  Finally, the correct constant for \textbf{(P2)} comes from the fact that
  $2A^\Phi=B^\Phi$ in that case.
\end{proof}

For any probability measure $\ga$ on $\dN$, and any $t\in\dR_+$, we denote by
$\ga\SGf{t}$ the probability measure on $\dN$ defined for any bounded function
$g:\dN\to\dR$ by $ \moy{\ga\SGf{t}}{g}:= \moy{\ga}{\SG{t}{g}}$. In particular,
when $\ga=\de_n$ for some fixed $n\in\dN$, we get
$\de_n\SGf{t}=\SG{t}{\cdot}(n)$. We have $\ga\ll\cP(\rho)$ for any probability
measure $\ga$ on $\dN$, as soon as $\rho>0$. Let us define
$f_\ga:=d\ga/d\cP(\rho)$. Since $\cP(\rho)$ is symmetric for $\GI$, we have
that $\GI$ and $\SG{t}{\cdot}$ are self-adjoint in $\rL^2(\cP(\rho))$.
Therefore, one can write for any $g\in\rL^2(\cP(\rho))$
$$
\int_\dN\!\SG{t}{f_\ga}g\,d\cP(\rho)
= \int_\dN\!\SG{t}{g}f_\ga\,d\cP(\rho)
= \int_\dN\!\SG{t}{g}\,d\ga
= \int_\dN\!g\,d(\ga\SGf{t}).
$$
Recall that the total variation norm $\NRM{\al}_\mathrm{TV}$ of a Borel
measure $\al$ on an at most countable set $\cS$ is defined by
$\NRM{\al}_\mathrm{TV} = \frac{1}{2}\NRM{\al}_1 =
\frac{1}{2}\sum_{x\in\cS}\ABS{\al(x)}$. If $\al$ and $\be$ are two probability
measures on $\cS$, the distance $\NRM{\al-\be}_\mathrm{TV}$ is
$$
\NRM{\al-\be}_\mathrm{TV} = \sup_{A\subset\cS}\ABS{\al(A)-\be(B)}
=\frac{1}{2}\sup_{\NRM{f}_\infty\leq1}\ABS{\int\!f\,d\al-\int\!f\,d\be}.
$$
Recall the well known bound for any $a,b\in\dR_+$,
$\NRM{\cP(a)-\cP(b)}_\mathrm{TV} \leq 1-e^{-(b-a)}$, cf. \cite[Prop. 6.1 page
143]{MR1996883}, which gives from \eqref{eq:mehler} for any $t\in\dR_+$
$$
\NRM{\SG{t}{\cdot}(0)-\cP(\rho)}_\mathrm{TV} %
\leq 1-e^{-\rho e^{-\mu t}}.
$$
Theorem \ref{th:MLSI} for \textbf{(P1)} produces in particular a bound for
$\NRM{\SG{t}{\cdot}(n)-\cP(\rho)}_\mathrm{TV}$, as stated in Corollary
\ref{co:ent-vt}. 

\begin{cor}\label{co:ent-vt}
  Let $(\SGf{t})_{t\geq0}$ be the semi-group of the M/M/$\infty$ queue
  with input rate $\la>0$ and service rate $\mu>0$.
  For any $n\in\dN$ and any $t\in\dR_+$,
  $$
  2\NRM{\SG{t}{\cdot}(n)-\cP(\rho)}_\mathrm{TV}^2
  \leq e^{-\mu t}\log(e^\rho\rho^{-n}n!).
  $$
\end{cor}

The proof follows the lines of a method due to Diaconis and Saloff-Coste, cf.
for example \cite{MR1410112,MR1490046}.

\begin{proof}[Proof of Corollary \ref{co:ent-vt}]\label{co:ent-vt:proof}
  Since $Q:=\cP(\rho)$ is symmetric for $\GI$,
  \begin{equation}\label{eq:nu-sg-radon-niko}
    \frac{d(\ga\SGf{t})}{dQ}=\SG{t}{f_\ga},
  \end{equation}
  where $f_\ga:=d\ga/dQ$. The Pinsker-Csisz\'ar-Kullback inequality states
  that for any couple $(\al,\be)$ of probability measures on the same measured
  space $2\NRM{\al-\be}_\mathrm{TV}^2 \leq\ent{}{\al\,\vert\,\be}
  =\ent{\be}{d\al/d\be}$, where $\entf{\be}$ is the $\Phi$-entropy in the case
  \textbf{(P1)}. Let $t\in\dR_+$ and $n\in\dN$. For
  $(\al,\be)=(\SG{t}{\cdot}(n),Q)$, we can write by
  \eqref{eq:nu-sg-radon-niko} and \eqref{eq:entropy-dec}
  $$
  2\NRM{\ga\SGf{t}-Q}_\mathrm{TV}^2
  \leq\ent{Q}{\SG{t}{f_\ga}}
  \leq e^{-\mu t}\ent{Q}{f_\ga}.
  $$
  For $\ga=\de_n$ for some fixed $n\in\dN$, we get
  $\ga\SGf{t}=\SG{t}{\cdot}(n)$ and $f_{\de_n}=\rI_{\{n\}}/Q(n)$. As a
  consequence, we obtain as expected $\ent{Q}{f_{\de_n}} = -\log
  Q(n)=\log(e^\rho\rho^{-n}n!)$.
\end{proof}

\subsection{Local inequalities and semi-group interpolation}

Standard Brownian motion on $\dR$ starting from $x$ interpolates on the time
interval $[0,t]$ between the Dirac measure $\de_x$ and the Gaussian measure
$\cN(x,t)$. It is known that this interpolation provides the optimal Gaussian
logarithmic Sobolev inequality. Similarly, the simple Poisson process of
intensity $\la$ starting from $n$ interpolates on the time interval $[0,t]$
between the Dirac measure $\de_n$ and the translated Poisson measure
$\de_n*\cP(\la t)$. By analogy, let us give a proof of the Poisson entropic
inequality \eqref{eq:pois-A} by using the simple Poisson process, which
corresponds to an M/M/$\infty$ queue with $\mu=0$. In that case, $\GI=\la\rD$
and $\SG{t}{\cdot}(0)=\cP(\la t)$. Let $\Phi\in\cC_\cI$ such that
$A^\Phi\in\cC_{\cT_\cI}$, cf. Theorem \ref{th:cvx-2}. One can write by
abridging $\SG{t}{\cdot}(0)$ in $\SG{t}{\cdot}$ and denoting $F=\SG{t-s}{f}$,
\begin{align*}
  \ENT{\cP(\la)}{\Phi}{f}
  & = \SG{1}{\Phi(f)}-\Phi(\SG{1}{f}) \\
  & = \int_0^t\!\pd_s\SG{s}{\Phi(\SGf{t-s}{f})}\,ds \\
  & = \int_0^t\!\SG{s}{\GI\Phi(F)-\Phi'(F)\GI F}\,ds.
\end{align*}
Now, $\GI\Phi(F)-\Phi'(F)\GI F=A^\Phi(F,\rD F)$, and \eqref{eq:commut} with
$\mu=0$ gives $\rD F=\rD\SG{t-s}{f}=\SG{t-s}{\rD f}$. Thus, we get,
\begin{align*}
  \ENT{\cP(\la)}{\Phi}{f}
  & = \la\int_0^t\SG{s}{A^\Phi(F,\rD F)}\,ds \\
  & = \la\int_0^t\SG{s}{A^\Phi(\SGf{t-s}{f},\SGf{t-s}{\rD f})}\,ds.
\end{align*}
Finally,
Jensen inequality for convex function $A^\Phi$ gives then the desired result,
$$
\ENT{\cP(\la)}{\Phi}{f} %
\leq \la\int_0^t\SG{s}{\SG{t-s}{A^\Phi(f,\rD f)}}\,ds %
= \la t\SG{t}{A^\Phi(f,\rD f)}.
$$

\subsubsection*{M/M/$\infty$ semi-group interpolation  on the time interval $[0,+\infty]$}

The standard Ornstein-Uhlenbeck process on $\dR$ starting from $x$
interpolates on the time interval $[0,+\infty]$ between the Dirac measure
$\de_x$ and the standard Gaussian measure $\cN(0,1)$. It is known that this
interpolation provides the optimal Gaussian logarithmic Sobolev inequality.
Similarly, the M/M/$\infty$ queue with intensities $(\la,\mu)$ starting from
$n$ interpolates on the time interval $[0,+\infty]$ between the Dirac measure
$\de_n$ and the Poisson measure $\cP(\rho)$ where $\rho=\la/\mu$. Notice that
when $\la=0$, this interpolation holds between $\de_n$ and $\de_0$ (pure death
process). By analogy, let us give a proof of the Poisson entropic inequality
\eqref{eq:pois-A} by using the M/M/$\infty$ queue. Let $(\SGf{t})_{t\geq0}$ be
the M/M/$\infty$ queue semi-group with input rate $\la$ and service rate
$\mu$. Let $\Phi\in\cC_\cI$ such that $B^\Phi\in\cC_{\cT_\cI}$, cf. Theorem
\ref{th:cvx-2}. We denote by $Q$ the Poisson measure $\cP(\rho)$. For any
$f\in\cK(\dN,\cI)$, we write
\begin{align*}
  \ENT{Q}{\Phi}{f}
  & = +\int_\dN \PAR{\Phi(\SGf{0}{f})-\Phi(\SGf{\infty}{f})}\,dQ \\
  & = -\int_\dN \int_0^{\infty} \pd_t\Phi(\SGf{t}{f})\,dt\,dQ \\
  & = -\int_0^{\infty}\int_\dN \Phi'(\SGf{t}{f})\GI\SGf{t}{f}\,dQ\,dt \\
  & = \la\int_0^{\infty}\int_\dN B^\Phi(\SGf{t}{f},\rD\SGf{t}{f})\,dQ\,dt,
\end{align*}
where we used \eqref{eq:propBPoi} for the last equality. Now, the commutation
\eqref{eq:commut} yields
\begin{align*}
  \ENT{Q}{\Phi}{f}
  &=\la\int_0^{\infty} \int_\dN
  B^\Phi(\SGf{t}{f},e^{-\mu t}\SGf{t}{\rD f})\,dQ\,dt.
\end{align*}
Jensen inequality for $B^\Phi$ and $\SG{t}{\cdot}$ followed by the
invariance of $Q$ give
\begin{align*}
  \ENT{Q}{\Phi}{f}
  &\leq\la\int_0^{\infty} \int_\dN B^\Phi(f,e^{-\mu t}\rD f)\,dQ\,dt.
\end{align*}
But by Lemma \ref{le:cvx-4}, $B^\Phi(u,e^{-\mu t}v)\leq e^{-\mu t}B^\Phi(u,v)$
for any $(u,v)\in\cT_\cI$, and thus
$$
\ENT{Q}{\Phi}{f}
\leq\la\int_0^{\infty}e^{-\mu t}\,dt \int_\dN B^\Phi(f,\rD f)\,dQ 
=\rho\moy{Q}{B^\Phi(f,\rD f)},
$$
which is exactly the $B$ transform version of the Poisson entropic inequality
\eqref{eq:pois-A}.

\begin{rem}[$A-B-C$ transforms and discrete space]
  The interpolation on $[0,t]$ gives rise to the $A^\Phi$ transform whereas
  the interpolation on $[0,+\infty]$ leads to the $B^\Phi$ transform. In
  continuous space settings, the diffusion property permits to write
  $\GI\Phi(F)-\Phi'(F)\GI F = \Phi''(F)\GA(F,F)$ which is close to $C^\Phi$
  and not to $A^\Phi$ in that case.
\end{rem}

\subsubsection*{Local inequality and semi-group interpolation on the time
  interval $[0,t]$}

Consider the semi-group $(\SGf{t})_{t\geq0}$ of the M/M/$\infty$ queue with
input rate $\la$ and service rate $\mu$. The family
$(\SGf{s}(\cdot)(n))_{0\leq s\leq t}$ interpolates between $\de_n$ and
$\cB(n,e^{-\mu t})*\cP(\rho(1-e^{-\mu t}))$. Let $\Phi\in\cC_\cI$ such that
$A^\Phi\in\cC_{\cT_\cI}$, cf. Theorem \ref{th:cvx-2}. The inequalities
\eqref{eq:binpois-A} and \eqref{eq:mehler} give for any $n\in\dN$, any
$t\in\dR_+$, and any $f\in\cK(\dN,\cI)$,
\begin{equation}\label{eq:mmi-loc}
  \ENT{\SG{t}{\cdot}(n)}{\Phi}{f}
  \leq \rho q(t)\SG{t}{A^\Phi(f,\rD f)}(n) %
  + np(t)q(t)\SG{t}{q(t)A^\Phi(f,\rD f)+p(t)A^\Phi(\tau(f,\rD f))}(n-1).
\end{equation}
Let us try to recover \eqref{eq:mmi-loc} by semi-group interpolation. We write
as for the pure Poisson process case,
\begin{align*}
  \ENT{\SGf{t}(\cdot)(n)}{\Phi}{f}
  & = \SG{t}{\Phi(f)}(n)-\Phi(\SG{t}{f}(n)) \\
  & = \int_0^t\!\pd_s\SG{s}{\Phi(\SGf{t-s}{f})}(n)\,ds \\
  & = \int_0^t\!\SG{s}{\GI\Phi(F)-\Phi'(F)\GI F}(n)\,ds.
\end{align*}
where $F:=\SG{t-s}{f}$. At this step, we notice that
$$
\GI\Phi(F)-\Phi'(F)\GI F=\la A^\Phi(F,\rD F)+\mu hA^\Phi(F,\rD^*F),
$$
where $h:\dN\to\dN$ is defined by $h(k)=k$ for any $k\in\dN$. Thus, we get
\begin{equation}\label{eq:ent-loc}
  \ENT{\SGf{t}(\cdot)(n)}{\Phi}{f} %
    = \la\int_0^t\!\SG{s}{A^\Phi(F,\rD F)}(n)\,ds %
    +\mu\int_0^t\!\SG{s}{hA^\Phi(F,\rD^* F)}(n)\,ds.
\end{equation}
By virtue of \eqref{eq:commut}, \eqref{eq:Ap-C-A}, Jensen inequality for the
convex functions $A^\Phi$ and $C^\Phi$, and the semi-group property,
the first term of the right hand side of \eqref{eq:ent-loc} is bounded above by
$$
\PAR{\frac{1}{2}\la\int_0^t\!p(t-s)^2q(t-s)\,ds}\SG{t}{C^\Phi(f,\rD f)}(n)%
+\PAR{\la\int_0^t\!p(t-s)^3\,ds}\SG{t}{A^\Phi(f,\rD f)}(n).
$$
For the second term of the right hand side of \eqref{eq:ent-loc}, we first
write by virtue of \eqref{eq:ipp-bin-poi-sg} and \eqref{eq:adtau},
$$
\mu\SG{s}{hA^\Phi(F,\rD^* F)}(n)
=\mu np(s)\SG{s}{A^\Phi(\tau(F,\rD F))}(n-1)
+\la q(s)\SG{s}{A^\Phi(\tau(F,\rD F))}(n).
$$
Now, by \eqref{eq:commut}, \eqref{eq:Atp-C-A}, Jensen inequality for the
convex functions $A^\Phi(\tau)$ and $C^\Phi$, and the semi-group property,
$$
\SG{s}{A^\Phi(\tau(F,\rD F))}(k)
\leq
\frac{1}{2}p(t-s)^2q(t-s)\SG{t}{C^\Phi(f,\rD f)}(k)
+p(t-s)^3\SG{t}{A^\Phi(\tau(f,\rD f))}(k)
$$
for any $k\in\{n-1,n\}$. Thus, the second term of the right hand side of
\eqref{eq:ent-loc} is bounded above by
\begin{multline*}
\PAR{\mu n\int_0^t\!p(s)p(t-s)^3\,ds}\SG{t}{A^\Phi(\tau(f,\rD f))}(n-1) 
+\PAR{\frac{1}{2}\mu %
n\int_0^t\!p(s)p(t-s)^2q(t-s)\,ds}\SG{t}{C^\Phi(f,\rD f)}(n-1) \\
+\PAR{\la\int_0^t\!q(s)p(t-s)^3\,ds}\SG{t}{A^\Phi(\tau(f,\rD f))}(n) 
+\PAR{\frac{1}{2}\la\int_0^t\!q(s)p(t-s)^2q(t-s)\,ds}\SG{t}{C^\Phi(f,\rD
  f)}(n).
\end{multline*}
Putting all together, we obtain finally the following local inequality.
\begin{align}\label{eq:mmi-loc-new}
  \ENT{\SGf{t}(\cdot)(n)}{\Phi}{f}
  \leq 
  \rho\SG{t}{\frac{1}{3}(1-p(t)^3)A^\Phi(f,\rD f) +\frac{1}{6}q(t)^2(2+p(t))
    \SBRA{A^\Phi(\tau(f,\rD f))+\frac{1}{2}C^\Phi(f,\rD f)}}(n)\nonumber\\
  +\frac{1}{2}np(t)\SG{t}{(1-p(t)^2)A^\Phi(\tau(f,\rD
    f))+\frac{1}{2}q(t)^2C^\Phi(f,\rD f)}(n-1),
\end{align}
which is not \eqref{eq:mmi-loc}. Actually, \eqref{eq:mmi-loc-new} is in a way
stronger than \eqref{eq:mmi-loc} for small $t$, as we will see in the sequel
with the fluid limit approximation of the Ornstein-Uhlenbeck process. When
$t\to\infty$, we have $p(t)\to0$, $q(t)\to1$, and
$A^\Phi+A^\Phi(\tau)=B^\Phi$, and in that case, \eqref{eq:mmi-loc-new}
provides the following Poissonian inequality.
$$
\ENT{\cP(\rho)}{\Phi}{f}
\leq \frac{1}{2}\rho\moy{\cP(\rho)}
{\frac{2}{3}B^\Phi(f,\rD f)+\frac{1}{3}C^\Phi(f,\rD f)}(n),
$$
which is not \eqref{eq:pois-A}. The proof of \eqref{eq:mmi-loc-new} given
above suggests to use \eqref{eq:AB-int-C} instead of its consequences
\eqref{eq:Ap-C-A} and \eqref{eq:Atp-C-A} for the derivation of local
inequalites via semi-group interpolation. The investigation of this approach
is left to the reader. Notice that \eqref{eq:pA-Ap} is not strong enough. Let
us focus on the \textbf{(P2)} case, for which we have the simple identity
$2A^\Phi(f,\rD f)=2A^\Phi(\tau(f,\rD f))=B^\Phi(f,\rD f)=C^\Phi(f,\rD f)
=2\ABS{\rD f}^2$. In that case, \eqref{eq:mmi-loc-new} is the optimal local
Poincar\'e inequality, e.g.
$$
\var{\SG{t}{\cdot}(n)}{f} 
\leq \rho q(t)\SG{t}{\ABS{\rD f}^2}(n)
+ n p(t)q(t)\SG{t}{\ABS{\rD f}^2}(n-1).
$$

\subsection{Scaling limit of the entropic inequalities}

Let us consider the Poisson distribution $\cP(\rho)$ with parameter $\rho>0$.
For any $N\in\dN^*$, let $\kappa_N:\dN\to\dR$ be the function defined by
$\kappa_N(n):=N^{-1/2}(n-\rho N)$ for any $n\in\dN$. By virtue of the Central
Limit Theorem, the image measure of $\cP(N\rho)=\cP(\rho)^{*N}$ by $\kappa_N$
converges weakly towards the Gaussian measure $\cN(0,\rho)$ when $N\to\infty$.
Let $g\in\cK(\dR,\cI)$ be smooth with bounded derivatives, and set
$f_N:=g\circ\kappa_N$. In one hand, we have
$$
\lim_{N\to\infty}\ENT{\cP(N\rho)}{\Phi}{f_N} = \ENT{\cN(0,\rho)}{\Phi}{g}.
$$
In the other hand, by a Taylor formula, $D(f_N)=\rD(g\circ\kappa_N) =
N^{-1/2}(g'\circ \kappa_N) + O(N)$, and by a Taylor formula for $\Phi$ this
time, $A^\Phi(f_N,\rD f_N) = (2N)^{-1}(g'\circ\kappa_N)^2\Phi''(f_N)+o(N)$.
This yields that
$$
\lim_{N\to\infty}\rho N\moy{\cP(\rho)}{A^\Phi(f_N,\rD f_N)}
=\frac{1}{2}\rho\moy{\cN(0,\rho)}{C^\Phi(g,g')}.
$$
Now, the $A$-transform based Poisson entropic inequality \eqref{eq:pois-A} for
$\cP(N\rho)$ and $f_N$ gives finally that
\begin{equation}\label{eq:phisob-opt-g}
\ENT{\cN(0,\rho)}{\Phi}{g} \leq \frac{1}{2}\rho\moy{\cN(0,1)}{C^\Phi(g,g')}.
\end{equation}
Recall that the Poincar\'e inequality corresponds to \textbf{(P2)}. In that
case,
$$
\ENT{\cN(0,\rho)}{\Phi}{g}=\var{\cN(0,\rho)}{g} \text{\quad and\quad}
C^\Phi(g,g')=2\ABS{g'}^2.
$$ 
The logarithmic Sobolev inequality corresponds to \textbf{(P1)}. In that case,
$$
\ENT{\cN(0,\rho)}{\Phi}{g}=\ent{\cN(0,\rho)}{g} \text{\quad and \quad}
C^\Phi(g,g')=\frac{\ABS{g'}^2}{g}.
$$ 
The constant $\rho$ in \eqref{eq:phisob-opt-g} is known to be optimal. It
gives in particular the optimal Poincar\'e inequality for the Gaussian measure
in the case \textbf{(P2)}, and the optimal logarithmic Sobolev inequality for
the Gaussian measure in the case \textbf{(P1)}. The method was used in the
case \textbf{(P1)} in \cite[Remark 1.6]{MR1800540}. In some sense, the $A$
transform is the right Dirichlet form to consider since it allows the
derivation of optimal Gaussian entropic inequalities from their $A$-transform
based Poisson versions. In contrast, it is shown in \cite[pages
356-357]{MR1636948} that the optimal $B$ transform version for the Poisson
measure does not lead to the optimal constant in the logarithmic Sobolev
inequality for the Gaussian measure (lack of a multiplicative factor $2$). The
deep reason for this difference between $A$ and $B$ transforms consequences is
the fact that the comparison $A^\Phi\leq B^\Phi$ improves by a factor $2$ when
$v$ goes to $0$, as stated in Remark \ref{rm:comp-A-B-C}. This phenomenon does
not hold for the Poincar\'e inequality, since $2A^\Phi=B^\Phi$ for
\textbf{(P2)}.

As presented in Section \ref{ss:ou-mmi-fluid}, the M/M/$\infty$ queueing
process gives rise to an Ornstein-Uhlenbeck process via a fluid limit
procedure. It is quite natural to ask about the behaviour of the
binomial-Poisson entropic inequalities under this scaling limit.

Let $(X_t^N)_{t\geq0}$ be an M/M/$\infty$ queueing process with input rate
$N\la>0$ and service rate $\mu>0$, where $N\in\dN^*$. Let $(U_t)_{t\geq0}$ be
an Ornstein-Uhlenbeck process, solution of the Stochastic Differential
Equation $dU_t=\la dB_t-\mu U_t\,dt$. Let $g\in\cK(\dR,\cI)$ be smooth with
bounded derivatives. For any $y\in\dR$, we define $z_N:=\PENT{N\rho+N^{1/2}y}$
where $\PENT{\cdot}$ denotes the integer part. According to Section
\ref{ss:ou-mmi-fluid}, the image measure of $\cL(X_t^N\,\vert\,X_0=z_N)$ by
function $\kappa_N$ converges weakly towards $\cL(U_t\,\vert\,U_0=y)$ when $N$
goes to $\infty$. Notice that
$$
\cL(X_t^N\,\vert\,X_0=z_N) = \cB(z_N,p(t))*\cP(N\rho q(t))
$$
and that $\cL(U_t\,\vert\,U_0=y) = \cN(yp(t),\rho(1-p(t)^2))$. In particular,
if $f_N:=g\circ\kappa_N$, then
$$
\lim_{N\to\infty}\ENT{\cL(X_t^N\,\vert\,X_0^N=z_N)}{\Phi}{f_N}%
=\ENT{\cL(U_t\,\vert\,U_0=y)}{\Phi}{g}.
$$
In the other hand, as for the pure Poisson measure case, we have
$$
\lim_{N\to\infty}N\moy{\cL(X_t^N\,\vert\,X_0^N=z_N)}{A^\Phi(f_N,\rD f_N)}
=\frac{1}{2}\moy{\cL(U_t\,\vert\,U_0=y)}{C^\Phi(g,g')}.
$$
A similarly identity holds for $A^\Phi(\tau(f_N,\rD f_N))$. Putting all
together, we deduce from \eqref{eq:mmi-loc} that
\begin{equation}\label{eq:phisob-g-mmi}
\ENT{\cL(U_t\,\vert\,U_0=y)}{\Phi}{g} 
\leq K(t)\moy{\cL(U_t\,\vert\,U_0=y)}{C^\Phi(g,g')},
\end{equation}
where $K(t):=\mbox{$\frac{1}{2}$}\rho q(t)(1+2p(t))$. It is known that the
best constant in \eqref{eq:phisob-g-mmi} is $K^*(t):=\mbox{$\frac{1}{2}$}\rho
q(t)(1+p(t))$. Let us consider now the function $\te:\dR_+\to\dR_+$ defined
for any $t\in(0,\infty)$ by
$$
\te(t):=\frac{K(t)}{K^*(t)}
=1+\frac{1}{1+\frac{1}{p(t)}}.
$$
This function is non-increasing, with $\te(0)=\mbox{$\frac{3}{2}$}$ and
$\lim_{t\to+\infty}\te(t)=1$. Consequently, the constant $K(t)$ in the
inequality \eqref{eq:phisob-g-mmi} improves when $t$ increases. It is
asymptotically optimal, when $t$ goes to $+\infty$. Surprisingly, it turns out
that the usage of \eqref{eq:mmi-loc-new} instead of \eqref{eq:mmi-loc}
provides \eqref{eq:phisob-g-mmi} with constant $K^*(t)$ instead of constant
$K(t)$. As a consequence, \eqref{eq:mmi-loc-new} is in a way stronger than
\eqref{eq:mmi-loc}, at least in terms of their fluid limit.

\begin{rem}[The M/M/$1$ case]
  The M/M/$1$ queue with input rate $\la$ and service rate $\mu$ is the birth
  and death process on $\dN$ with generator $\GI=\mu\rD^*+\la\rD$. We have
  $[\GI,\rD]=0$ and the ``curvature'' is identically zero. When $\la>0$ and
  $\mu>0$, the symmetric invariant measure $Q$ is given by $Q(n)=\rho^n$ for
  any $n\in\dN$, with $\rho:=\la/\mu$. The associated Markov semi-group
  $(\SGf{t})_{t\geq0}$ satisfies to the exact commutation formula
  $\rD\SGf{t}=\SG{t}{\rD}$. Measure $Q$ is finite if and only if $\rho\leq 1$,
  and $Q$ is in that case the geometric measure $\cG(1-\rho)$ of mean
  $\la/(\mu-\la)$. This leads to Poisson like entropic inequalities for
  $\SG{t}{\cdot}$. The M/M/$1$ is a discrete space analog of the continuous
  process $(E_t)_{t\geq0}$ solution of the Stochastic Differential Equation
  $dE_t=dB_t-\mathrm{sign}(E_t)dt$. 
\end{rem}

\begin{rem}[Spectrum]
  A function $f:\dN\to\dR$ is an eigenvector associated to the eigenvalue
  $\al\in\dR$ for the M/M/$\infty$ infinitesimal generator $\GI$ defined by
  \eqref{eq:mmi-gi} if and only if $\la f(n+1)=(\la+\al+n\mu)f(n)-n\mu f(n-1)$
  for any $n\in\dN$. Obviously, for any $\al\in\dR$ and any starting value
  $f(0)\neq0$, the equation above has a unique non null solution. As a
  consequence, the spectrum of $\GI$ is $\dR$. We will denote by $f_\al$ the
  unique solution such that $f_\al(0)=1$. By the equation $\GI f_\al = \al f$
  and the invariance of $Q:=\cP(\rho)$ we get that $\moy{Q}{f_\al}=0$ as soon
  as $f_\al\in\rL^1(Q)$. Suppose that $f_\al\in\rL^2(Q)$, then
  $0\leq\moy{Q}{\GA(f,f)}=-\moy{Q}{f\GI f}=-\al\moy{Q}{f^2}$ and thus
  $\al\leq0$. Moreover, Theorem \ref{th:MLSI} for \textbf{(P2)} gives
  $\inf\BRA{-\al\in\dR,\,f_\al\in\rL^2(Q)}=\mu^{-1}$, cf. \cite[Prop.
  2.3]{MR1307413}.
\end{rem}

\begin{rem}[Bakry $\GD$ calculus]
  Let us define the Markovian functional quadratic forms $\GA$ and $\GD$ by
  $2\GA(f,f):=\GI(f^2)-2f\GI f$ and $2\GD(f,f):=\GI\GA(f,f)-2\GA(f,\GI f)$.
  After some algebra based on \eqref{eq:mmi-gi}, we get
  \begin{equation*}
    2\GA(f,f)(n)=n\mu\ABS{\rD^* f}^2(n)+\la\ABS{\rD f}^2(n)
  \end{equation*}
  for any $f:\dN\to\dR$ and any $n\in\dN$, and
  \begin{equation*}
    2\GD(f,f)(n)=\frac{3}{2}\la\mu\ABS{\rD f}^2(n)
    +\frac{n}{2}\mu^2\ABS{\rD^*f}^2(n) +\bR(f,f)(n),
  \end{equation*}
  where $2\bR(f,f)(n):= n(n-1)\mu^2\ABS{\rD^*\rD^*f}^2
  +2n\la\mu\ABS{\rD\rD^*f}^2 +\la^2\ABS{\rD\rD f}^2$. Notice that for the
  linear function $f(n)=n$, we get
  $$
  2\GA(f,f)(n)=\la+n\mu\text{\quad and\quad} 4\GD(f,f)(n)=3\la\mu+n\mu^2.
  $$
  Since $\bR(f,f)\geq0$ for any $f$, we obtain immediately the bound
  $\GD\geq\mu\frac{1}{2}\GA$, which is the infinitesimal version of the
  commutation $\GA\SGf{t}\leq \exp(-t\frac{\mu}{2})\SGf{t}{\GA}$. Moreover, an
  integration by parts similar to \eqref{eq:ipp-poi} gives the integrated
  bound $\moy{Q}{\GD f}\geq\mu\moy{Q}{\GA f}$, where $Q:=\cP(\rho)$. Such a
  bound gives, via integration by parts, the Poincar\'e inequality
  $\var{Q}{f}\leq\rho\moyf{Q}(\ABS{\rD f}^2)$, which is exactly
  \eqref{eq:pois-A} for \textbf{(P2)}. However, the $\GD$ bound above suggests
  that $\GD$ is not the right tool in order to derive $\Phi$-entropic
  inequalities beyond the \textbf{(P2)} case. Bakry-\'Emery type approaches
  are designed for diffusion. In discrete space settings, the lack of chain
  rule reduces their strength for the derivation of entropic inequalities
  beyond the \textbf{(P2)} case.
\end{rem}

\section{Convexity and $\Phi$-calculus on $A-B-C$ transforms }
\label{se:cvx}

We give in the sequel various convexity properties, which extend in particular
many aspects of \cite{MR2081075}. 
Let $\Phi:\cI\to\dR$ be a smooth function defined on an open interval
$\cI\subset\dR$. The usage of suitable Taylor formulas provide for any
$(u,v)\in\cT_\cI$,
\begin{equation}\label{eq:AB-int-C}
  A^\Phi(u,v)=\int_0^1\!(1-p)C^\Phi(u+pv,v)\,dp
  \text{\quad and \quad}
  B^\Phi(u,v)=\int_0^1\!C^\Phi(u+pv,v)\,dp,
\end{equation}
and for any $(u,v)\in\cT'_\cI$ and small enough $\veps$,
\begin{equation}\label{eq:AB-dif-C}
  A^\Phi(u,\veps v)=\frac{1}{2}C^\Phi(u,v)\veps^2+o(\veps^2)
  \text{\quad and \quad}
  B^\Phi(u,\veps v)=C^\Phi(u,v)\veps^2+o(\veps^2).
\end{equation}
We denote by $\tau:\cT_\cI\to\cT_\cI$ the bijective linear map defined for any
$(u,v)\in\cT_\cI$ by 
\begin{equation}\label{eq:def-tau}
\tau(u,v):=(u+v,-v).
\end{equation} 
Notice that $\tau$ is well defined since $(u,v)\in\cT_\cI$ implies that
$(u+v,u+v-v)\in\cI\times\cI$ and thus $(u+v,-v)\in\cT_\cI$. Writing
$(a,b):=(u,u+v)$ shows that the map $\tau$ transposes $a$ and $b$, and
$\tau^2$ is the identity map. Moreover,
\begin{align*}
 A^\Phi(u,v) & =\Phi(b)-\Phi(a)-\Phi'(a)(b-a), \\
 B^\Phi(u,v) & =(b-a)(\Phi'(b)-\Phi'(a)), \\
 C^\Phi(u,v) & =\Phi''(a)(b-a)^2.
\end{align*}

\begin{lem}\label{le:cvx-1}
  Let $\Phi:\cI\to\dR$ be smooth on an open interval $\cI\subset\dR$. Then the
  following statements hold.
  \begin{enumerate}
  \item $A^\Phi+A^\Phi(\tau)=B^\Phi$ and
    $B^\Phi(\tau)=B^\Phi$;
  \item Each of $A^\Phi$, $B^\Phi$, $C^\Phi$ is
    non-negative if and only if $\Phi\in\cC_\cI$.
  \end{enumerate}
\end{lem}

\begin{proof}
  The first statement of the Lemma is immediate. For the second statement, we
  observe first that $C^\Phi$ is non-negative if and only if $\Phi''$ is
  non-negative, e.g if and only if $\Phi\in\cC_\cI$. The same holds then for
  $A^\Phi$ and $B^\Phi$ by using \eqref{eq:AB-dif-C} and \eqref{eq:AB-int-C}.
\end{proof}

Lemma \ref{le:cvx-1} tells that the $A-B-C$ transforms map the set of convex
functions on $\cI$ into the set of non-negative functions on $\cT_\cI$.
Moreover, their null space contains any real valued affine functions on $\cI$.

\begin{lem}\label{le:cvx-3/2}
  Let $\Phi:\cI\to\dR$ be smooth on an open interval $\cI\subset\dR$. The
  following statements hold.
  \begin{enumerate}
  \item\label{le:cvx-3/2-1} for \textbf{(P1-P2-P3)},
    we have $\Phi''>0$ on $\cI$ and $\Phi$, $-\Phi'$, $\Phi''$, $-1/\Phi''$
    belong to $\cC_\cI$;
  \item\label{le:cvx-3/2-2} $2A^\Phi=B^\Phi=C^\Phi$ for \textbf{(P2)} and
    $A^\Phi\leq C^\Phi$ for \textbf{(P1)};
  \item\label{le:cvx-3/2-3} if $\Phi\in\cC_\cI$ then $A^\Phi\leq B^\Phi$.
  \item\label{le:cvx-3/2-4} if $\Phi''\in\cC_\cI$ then $C^\Phi(u+v/3,v)\leq
    2A^\Phi(u,v)$ and $C^\Phi(u+v/2,v)\leq B^\Phi(u,v)$ for any
    $(u,v)\in\cT_\cI$.
\end{enumerate}
\end{lem}

\begin{proof}
  Statement \ref{le:cvx-3/2-1} and the first part of statement
  \ref{le:cvx-3/2-2} are immediate. Notice that $1/\Phi''$ is affine for
  \textbf{(P1)} and \textbf{(P2)}. The second part of statement
  \ref{le:cvx-3/2-2} follows from the first part of \eqref{eq:AB-int-C}. For
  statement \ref{le:cvx-3/2-3}, we notice that by Lemma \ref{le:cvx-1},
  $B^\Phi=A^\Phi+A^\Phi(\tau)$, where $A^\Phi(\tau)\geq0$ when
  $\Phi\in\cC_\cI$. Statement \ref{le:cvx-3/2-4} follows by using
  \eqref{eq:AB-int-C}, the definition of $C^\Phi$, and Jensen inequality with
  respect to the integral over $[0,1]$ for the convex function
  $p\in[0,1]\mapsto\Phi''(u+pv)$.
\end{proof}

\begin{rem}[Optimality of $A$-$B$-$C$ comparisons]\label{rm:comp-A-B-C}
  The bound $A^\Phi\leq B^\Phi$ is optimal in the sense that for
  \textbf{(P1)}, we have $B^\Phi(u,v)\sim A^\Phi(u,v)$ at $v=+\infty$ for any
  $u\in\cI$. However, $B^\Phi=2A^\Phi=C^\Phi$ for \textbf{(P2)}; and in
  general
  $$
    \lim_{v\to0}v^{-2}2A^\Phi(u,v) %
    =\lim_{v\to0}v^{-2}B^\Phi(u,v) %
    =v^{-2}C^\Phi(u,v) %
    =\Phi''(u).
  $$
\end{rem}

Theorem \ref{th:cvx-2} below states that the convexity of the $A-B-C$
transforms of $\Phi$ are deeply related to the convexity of the $\Phi$-entropy
functional. It provides in particular a synthesis of some results by Lata{\l}a
and Oleszkiewicz in \cite{MR1796718}, by the author in \cite{MR2081075}, and
by Massart in his Saint-Flour course \cite{massart-st-flour} (see also the
article \cite{MR2123200}). We say that a collection $\cP$ of probability
spaces is a \emph{covering collection} if
$\{Q(T);T\in\cE,(E,\cE,Q)\in\cP\}=[0,1]$. An example is given for instance by
the family of Bernoulli probability measures on the two point space $\{0,1\}$,
or by any collection containing a probability measure on $\dR$ with a
continuous cumulative distribution function.

\begin{thm}\label{th:cvx-2} 
  For any smooth $\Phi:\cI\to\dR$ on an open interval $\cI\subset\dR$, the
  following statements are equivalent.
  \begin{enumerate}
  \item\label{th:cvx-2:1} $A^\Phi\in\cC_{\cT_\cI}$;
  \item\label{th:cvx-2:2} $B^\Phi\in\cC_{\cT_\cI}$;
  \item\label{th:cvx-2:3} $C^\Phi\in\cC_{\cT'_\cI}$;
  \item\label{th:cvx-2:4} either $\Phi$ is affine on $\cI$, or $\Phi''>0$ on
    $\cI$ with $-1/\Phi''\in\cC_\cI$;
  \item\label{th:cvx-2:5} $(a,b)\in\cI\times\cI\mapsto
    t\Phi(a)+(1-t)\Phi(b)-\Phi(ta+(1-t)b)$ belongs to $\cC_{\cI\times\cI}$,
    for any $t\in[0,1]$;
  \item\label{th:cvx-2:6} for any probability space $(E,\cE,Q)$,
    $\ENTF{Q}{\Phi}\in\cC_{\cK(E,\cI)}$;
  \item\label{th:cvx-2:7} there exists a covering collection $\cP$ such that
    $\ENTF{Q}{\Phi}\in\cC_{\cK(E,\cI)}$ for any $(E,\cE,Q)$ in
    $\cP$;
  \item\label{th:cvx-2:8} for any probability space $(E,\cE,Q)$ and any
    $f\in\cK(E,\cI)$,
    \begin{equation}\label{eq:varfor}
      \ENT{Q}{\Phi}{f}= \sup_{g\in\cK(E,\cI)}
      \{\moy{Q}{(\Phi'(g)-\Phi'(\moyf{Q}g))(f-g)} + \ENT{Q}{\Phi}{g}\};
    \end{equation}
  \item\label{th:cvx-2:9} there exists a covering collection $\cP$ such that
    \eqref{eq:varfor} holds for any $(E,\cE,Q)\in\cP$ and any
    $f\in\cK(E,\cI)$;
  \item\label{th:cvx-2:10} for any product probability space
    $(E,\cE,Q):=(E_1\times\cdots\times
    E_n,\cE_1\otimes\cdots\otimes\cE_n,Q_1\otimes\cdots\otimes Q_n)$, and for
    any $f\in\cK(E,\cI)$,
    \begin{equation}\label{eq:tenso}
      \ENT{Q}{\Phi}{f}
      \leq \moy{Q}{\ENT{Q_1}{\Phi}{f}}+\cdots+\moy{Q}{\ENT{Q_n}{\Phi}{f}},
    \end{equation}
    where the expectation with respect to $Q_i$ in $\ENT{Q_i}{\Phi}{f}$
    concerns only the $i^\text{th}$ coordinate;
  \item\label{th:cvx-2:11} there exists a covering collection $\cP$ such that
    \eqref{eq:tenso} holds for $n=2$, any $Q_1\in\{\cB(1,p);p\in[0,1]\}$, any
    $Q_2\in\cP$, and any $f\in\cK(E_1\times E_2,\cI)$.
  \end{enumerate}
  Moreover, if these statements hold, then $\Phi$ and $\Phi''$ belong to
  $\cC_\cI$.
\end{thm}

\begin{rem}[Functional spaces]\label{rm:L1Phi}
  By approximation, the convex set $\cK(E,\cI)$ can be replaced by the convex
  set $\rL^{1,\Phi}(Q)$ in statements \ref{th:cvx-2:6}, \ref{th:cvx-2:7},
  \ref{th:cvx-2:8}, \ref{th:cvx-2:9}, \ref{th:cvx-2:10}, \ref{th:cvx-2:11} of
  Theorem \ref{th:cvx-2}. More precisely, statement \ref{th:cvx-2:4} implies
  the convexity of $\Phi$, which implies in turn that
  $\Phi'(g)(f-g)+\Phi(g)\leq\Phi(f)$ for any $f,g\in\rL^{1,\Phi}(Q)$. This
  yields that $\moy{Q}{(\Phi'(g)-\moyf{Q}{g})(f-g)}$ is well defined in
  $[-\infty,+\infty)$, as noticed in the proof of \cite[Lem.
  2.26]{massart-st-flour}.
\end{rem}

\begin{rem}[Meaning of the variational formula]\label{rm:varfor}
  Despite its functional expression, the variational formula \eqref{eq:varfor}
  is actually a unidimensional statement, taken in all directions. Namely, for
  any probability space $(E,\cE,Q)$ and any $f,g\in\rL^{1,\Phi}(Q)$, let us
  denote by $\al_{f,g}:[0,1]\to\dR$ the function defined for any $\la\in[0,1]$
  by
  \begin{equation}\label{eq:def-al-fg}
    \al_{f,g}(\la):=\ENT{Q}{\Phi}{\la f+(1-\la)g}.
  \end{equation}
  Notice that $\al_{f,g}(0)=\ENT{Q}{\Phi}{g}$ and
  $\al_{f,g}(1)=\ENT{Q}{\Phi}{f}$. The consideration of convex combinations
  reveals that the convexity of the $\Phi$-entropy functional on
  $\rL^{1,\Phi}(Q)$ is equivalent to the convexity of $\al_{f,g}$ on $[0,1]$
  for any $f$ and $g$. Assume now that $\al_{f,g}$ is convex on $[0,1]$ for
  any $f$ and $g$ in $\rL^{1,\Phi}(Q)$. Assume for the moment that
  $f,g\in\cK(E,\cI)$. In that case, there are no boundary effects, and
  $\al_{f,g}$ is smooth. Recall that a real convex function is the envelope of
  its tangents, cf. \cite{MR1451876}. In particular,
  $\al_{f,g}(1)\geq\al_{f,g}(0)+\al_{f,g}'(0)$. Moreover, equality is achieved
  for $f=g$. As a consequence, we get
  \begin{equation}\label{eq:varfor-al-fg}
    \ENT{Q}{\Phi}{f}%
    =\sup_{g\in\cK(E,\cI)}\BRA{\al'_{f,g}(0)+\al_{f,g}(0)}.
  \end{equation}
  It turns out that
  $\al_{f,g}'(0)=\moy{Q}{(\Phi'(g)-\Phi'(\moyf{Q}{g}))(f-g)}$. We thus recover
  exactly \eqref{eq:varfor}. By virtue of Remark \ref{rm:L1Phi}, the formula
  above for $\al_{f,g}'(0)$ still makes sense in $[-\infty,\infty)$ when $f,g$
  are in $\rL^{1,\Phi}(Q)$, and consequently, the variational formula
  \eqref{eq:varfor-al-fg} remains true when $f,g\in\rL^{1,\Phi}(Q)$. Notice
  that $\al_{f,g}(\la)=\ENT{Q}{\Phi}{g+\la(f-g)}$, and hence $\al'_{f,g}(0)$
  is the directional derivative of the $\Phi$-entropy functional at point $g$
  in the direction $f-g$.
\end{rem}

\begin{proof}[Proof of Theorem \ref{th:cvx-2}]
  \emph{\ref{th:cvx-2:1}$\Rightarrow$\ref{th:cvx-2:2}}. Follows from the
  identity $B^\Phi=A^\Phi+A^\Phi(\tau)$ where $\tau$ is linear, given by Lemma
  \ref{le:cvx-1}.

  \emph{\ref{th:cvx-2:3}$\Rightarrow$\ref{th:cvx-2:1}} and
  \emph{\ref{th:cvx-2:3}$\Rightarrow$\ref{th:cvx-2:2}}. Follow from
  \eqref{eq:AB-int-C} used on a convex combination.
  
  \emph{\ref{th:cvx-2:1}$\Rightarrow$\ref{th:cvx-2:3}} and
  \emph{\ref{th:cvx-2:3}$\Rightarrow$\ref{th:cvx-2:2}}. Follow from
  \eqref{eq:AB-dif-C} used on a convex combination.
 
  \emph{\ref{th:cvx-2:1}$\Rightarrow$\ref{th:cvx-2:4}}. The Hessian matrix of
  $A^\Phi$ writes for any $(u,v)\in\cT_\cI$,
  \begin{equation*}
    \GR^2A^\Phi(u,v)=
    \left(
      \begin{array}{rl}
        A^{\Phi''}(u,v) & \Phi''(u+v)-\Phi''(u) \\
        \Phi''(u+v)-\Phi''(u) & \Phi''(u+v)
      \end{array}
    \right).
  \end{equation*}
  Since $A^\Phi$ is convex, the diagonal elements of $\GR^2A^\Phi$ are
  non-negative, and thus $\Phi''\geq0$ on $\cI$. Moreover, the convexity of
  $A^\Phi$ yields that $\det(\GR^2A^\Phi)$ is non-negative. Suppose now that
  $(u,v)\in\cT_\cI$ is such that $\Phi''(u+v)=0$. Then $\det(\GR^2
  A^\Phi(u,v))=-\Phi''(u)^2$, and thus $\Phi''(u)=0$. Consequently, the set
  $\{w\in\cI;\Phi''(w)=0\}$ is either empty of equal to $\cI$, as required.
  When $\Phi''>0$ on $\cI$, we get $\det(\GR^2A^\Phi(u,v)) =
  \Phi''(u+v)\Phi''^2(u)A^{-1/\Phi''}(u,v)$, which is non-negative since
  $A^\Phi\in\cC_{\cT_\cI}$. But $\Phi''>0$, and thus $A^{-1/\Phi''}\geq0$.
  Lemma \ref{le:cvx-1} gives then $-1/\Phi''\in\cC_\cI$.

  \emph{\ref{th:cvx-2:4}$\Rightarrow$\ref{th:cvx-2:1}}. If $\Phi$ is affine on
  $\cI$, then $A^\Phi$ is identically zero, and thus belongs to
  $\cC_{\cT_\cI}$. Let us consider the second case. Assume that $\Phi''>0$ on
  $\cI$ with $-1/\Phi''\in\cC_\cI$. It turns out that
  $(-1/\Phi'')''=(\Phi''''\Phi''-2\Phi'''^2)/\Phi''^3$. Hence, $\Phi''''\Phi''
  \geq 2\Phi'''^2$ on $\cI$, and thus $\Phi''''\geq0$ on $\cI$. In other
  words, $\Phi''\in\cC_\cI$. By Lemma \ref{le:cvx-1}, it follows that
  $A^{\Phi''}$ is non-negative. Therefore, the diagonal elements of $\GR^2
  A^\Phi$ are non-negative on $\cC_\cI$. In the other hand, for any
  $(u,v)\in\cT_\cI$, $\det(\GR^2A^\Phi(u,v)) =
  \Phi''(u+v)\Phi''^2(u)A^{-1/\Phi''}(u,v)$, which is non-negative again by
  virtue of Lemma \ref{le:cvx-1}. Putting all together, the two dimensional
  matrix $\GR^2 A^\Phi(u,v)$ has a non-negative trace and determinant for any
  $(u,v)\in\cC_\cI$, as expected.

  \emph{\ref{th:cvx-2:4}$\Rightarrow$\ref{th:cvx-2:5}} and
  \emph{\ref{th:cvx-2:5}$\Rightarrow$\ref{th:cvx-2:6}}. Follow from the
  definitions. See for instance \cite{MR1796718}.
 
  \emph{\ref{th:cvx-2:6}$\Rightarrow$\ref{th:cvx-2:7}} and
  \emph{\ref{th:cvx-2:8}$\Rightarrow$\ref{th:cvx-2:9}} and
  \emph{\ref{th:cvx-2:10}$\Rightarrow$\ref{th:cvx-2:11}} are immediate.

  \emph{\ref{th:cvx-2:1}$\Rightarrow$\ref{th:cvx-2:8}}. Let $f$ and $g$ be in
  $\cK(E,\cI)$. Since $A^\Phi\in\cC_{\cT_\cI}$, the following sort of
  ``$A^\Phi$-entropy''
  \begin{equation}\label{eq:def-J}
    J(f,g):=\moy{Q}{A^\Phi(g,f-g)}-A^\Phi(\moyf{Q}{g},\moy{Q}{f-g})
  \end{equation}
  is non negative by Jensen inequality. Moreover, it vanishes when $f=g$. The
  desired result follows from the identity $\ENT{Q}{\Phi}{f} %
  = J(f,g) %
  +\moy{Q}{(\Phi'(g)-\Phi'(\moyf{Q}{g}))(f-g)} + \ENT{Q}{\Phi}{g}$.

  \emph{\ref{th:cvx-2:9}$\Rightarrow$\ref{th:cvx-2:1}}. For any
  $(E,\cE,Q)\in\cP$ and $f,g\in\cK(E,\cI)$, the identity \eqref{eq:varfor}
  implies that the quantity $J(f,g)$ defined by \eqref{eq:def-J} is
  non-negative. By approximation, $J(f,g)$ is non-negative for any
  $f,g\in\rL^{1,\Phi}(Q)$. Now, let $(u,v)$ and $(u',v')$ be in $\cT_\cI$, and
  let $\la\in[0,1]$. Since $\cP$ is a covering collection, there exists
  $(E,\cE,Q)\in\cP$ and $T\in\cE$ such that $Q(T)=\la$. Let
  $g:=u\rI_T+u'\rI_{T^c}$ and $f:=g+v\rI_T+v'\rI_{T^c}$. The desired result
  follows from the identity $J(f,g)=\la
  A^\Phi(u,v)+(1-\la)A^\Phi(u',v')-A^\Phi(\la(u,v)+(1-\la)(u',v'))$.

  \emph{\ref{th:cvx-2:3}$\Rightarrow$\ref{th:cvx-2:6}} and
  \emph{\ref{th:cvx-2:3}$\Rightarrow$\ref{th:cvx-2:8}}. Let
  $f,g\in\cK(E,\cI)$, and let $\al_{f,g}$ be as in \eqref{eq:def-al-fg}. It
  turns out that
  $$
  \al''(t) %
  = \moy{Q}{C^\Phi(h_t,f-g)}-C^\Phi(\moyf{Q}{h_t},\moy{Q}{f-g}).
  $$
  Notice that $(h_t,f-g)$ takes its values in $\cT'_\cI$. Since
  $C^\Phi\in\cC_{\cT'_\cI}$, we get $\al''(t)\geq 0$. In other words,
  $\al\in\cC_{[0,1]}$. In particular $\al(\la) \leq \la\al(1)+(1-\la)\al(0)$
  for $\la\in[0,1]$ writes $\ENT{Q}{\Phi}{\la f+(1-\la)g} %
  \leq \la\ENT{Q}{\Phi}{f} + (1-\la)\ENT{Q}{\Phi}{g}$, which is nothing else
  but the expression of the convexity of $\ENTF{Q}{\Phi}$. Additionally, since
  every convex function on an interval is the envelope of its tangents, see
  \cite{MR1451876}, one gets
  $\ENT{Q}{\Phi}{f}=\al(1)=\sup_{t\in[0,1]}\BRA{\al(t)+\al'(t)(1-t)}$. In
  particular, $\ENT{Q}{\Phi}{f}\geq\al(0)+\al'(0)$, with equality when $f=g$.
  Taking the supremum with respect to $g$ leads to \eqref{eq:varfor}.

  \emph{\ref{th:cvx-2:7}$\Rightarrow$\ref{th:cvx-2:9}} and
  \emph{\ref{th:cvx-2:6}$\Rightarrow$\ref{th:cvx-2:8}}. Let
  $f,g\in\cK(E,\cI)$, and let $\al_{f,g}$ be as in
  \eqref{eq:def-al-fg}. Then, for any $s,t\in[0,1]$ and any $\la\in[0,1]$,
  $\al(\la s+(1-\la)t) = \ENT{Q}{\Phi}{\la(tf+(1-t)g)+(1-\la)(sf+(1-s)g)}$.
  Since $\ENTF{Q}{\Phi}\in\cC_{\cK(E,\cI)}$, we get $\al(\la
  s+(1-\la)t)\leq\la\al(s)+(1-\la)\al(t)$, and thus $\al\in\cC_{[0,1]}$. Since
  every convex function on an interval is the envelope of its tangents, see
  \cite{MR1451876}, we obtain
  $\ENT{Q}{\Phi}{f}=\al(1)=\sup_{t\in[0,1]}\BRA{\al(t)+\al'(t)(1-t)}$. In
  particular, $\ENT{Q}{\Phi}{f}\geq\al(0)+\al'(0)$, with equality when $f=g$.
  Taking the supremum with respect to $g$ leads to \eqref{eq:varfor}.

  \emph{\ref{th:cvx-2:9}$\Rightarrow$\ref{th:cvx-2:7}} and
  \emph{\ref{th:cvx-2:8}$\Rightarrow$\ref{th:cvx-2:6}}. Use $\la
  f_1+(1-\la)f_2-g=\la(f_1-g)+(1-\la)(f_2-g)$ and
  $\ENT{Q}{\Phi}{g}=\la\ENT{Q}{\Phi}{g}+(1-\la)\ENT{Q}{\Phi}{g}$ in the
  expression inside the supremum in \eqref{eq:varfor}, then use the fact that
  the supremum of the sum is less than or equal to the sum of the suprema.
 
  \emph{\ref{th:cvx-2:11}$\Rightarrow$\ref{th:cvx-2:7}}. The proof can be
  found in \cite[introduction of section 2.5]{massart-st-flour}. Namely, let
  $g_1,g_2\in\cK(E_2,\cI)$, and consider $f:\{0,1\}\times E_2\to\dR$ defined
  by $f(x,y):=g_1(y)$ if $x=0$ and $f(x,y):=g_2(y)$ if $x=1$. The
  tensorisation formula \eqref{eq:tenso} expressed for $f$ rewrites
  $\ENT{Q_2}{\Phi}{(1-p)g_1+pg_2}%
  \leq(1-p)\ENT{Q_2}{\Phi}{g_1}+p\ENT{Q_2}{\Phi}{g_2}$, as expected.

  \emph{\ref{th:cvx-2:8}$\Rightarrow$\ref{th:cvx-2:10}}. The proof can be
  found in the Saint-Flour course \cite[Theorem 2.27]{massart-st-flour}.
  Roughly speaking, it consists in the usage of the variational formula
  \eqref{eq:varfor} on each entropy in the right hand side of
  \eqref{eq:tenso}, which gives rise, via a telescopic sum, to the variational
  formula for the left hand side of \eqref{eq:tenso}.

  Finally, if the statements hold, then $\Phi\in\cC_\cI$ by statement
  \ref{th:cvx-2:4}, and the proof of
  \emph{\ref{th:cvx-2:4}$\Rightarrow$\ref{th:cvx-2:1}} given above provides in
  particular that $\Phi''\in\cC_\cI$.
\end{proof}

\begin{xpl} 
  For \textbf{(P1-P2-P3)}, both $\Phi$, $-\Phi'$, and $\Phi''$ are convex on
  $\cI$. Moreover, $\Phi''>0$ on $\cI$ and $-1/\Phi''$ is convex on $\cI$.
  Actually, $-1/\Phi''$ is affine for \textbf{(P1)} and \textbf{(P2)}.
  Consider the case where $\Phi(u):=-u\log(-u)$ on $\cI=(-\infty,0)$. Then
  $\Phi''>0$ on $\cI$ , and $-1/\Phi''>0$ on $\cI$. However, $-\Phi'$ is
  concave and not convex on $\cI$. Consider now the case where
  $\Phi(u)=-\log(u)$ on $\cI=(0,+\infty)$. Then $\Phi$, $-\Phi'$, $\Phi''$ are
  convex on $\cI$, and $\Phi''>0$ on $\cI$. However, $-1/\Phi''$ is concave
  and not convex. These examples rely on the stability by symmetry of the
  convexity of $-1/\Phi''$, and the absence of such a stability for $-\Phi'$.
\end{xpl}

\begin{xpl}
  Following \cite{MR2081075}, the convexities of the $A$-$B$-$C$ transforms of
  $\Phi$ and of the $\Phi$-entropy functional are stable by any linear
  combination on $\Phi$ with non-negative coefficients. Theorem \ref{th:cvx-2}
  shows in particular that this stability still holds for the convexity of
  $-1/\Phi''$, for which it is less apparent. The consideration of continuous
  linear combinations on $\Phi$ by mean of an integral with respect to a
  positive Borel measure provides several interesting examples. For instance,
  $\Phi(u):=\int_1^2\!u^p\,dp=u(u-1)/\log(u)$ on $\cI=\dR_+^*$ is obtained
  from \textbf{(P3)}, and satisfies to the required convexities of Theorem
  \ref{th:cvx-2}.
\end{xpl}

\begin{xpl}
  A curious example is given by $\Phi(u)=-\bI(u)$ on $\cI=(0,1)$, where $\bI$
  is the Gaussian isoperimetric function defined by $\bI:=F'\circ F^{-1}$
  where $F$ is the cumulative distribution function of $\cN(0,1)$. Function
  $\bI$ is positive and concave on $\cI$, and satisfies to the identity
  $\bI\bI''=-1$. Consequently, $\Phi$ and $-1/\Phi''=\Phi$ are convex. In
  particular, Theorem \ref{th:cvx-2} shows that the $\bI$-entropy is concave,
  and provides a reversed tensorisation formula.
\end{xpl}

For any $p\in[0,1]$, let $\si_p:\cT_\cI\to\cT_\cI$ be the linear map defined
for any $(u,v)\in\cT_\cI$ by 
\begin{equation}\label{eq:def-sigp}
  \si_p(u,v):=(u,pv).
\end{equation}
The map $\si_p$ is well defined since for any $(u,v)\in\cT_\cI$ and any
$p\in[0,1]$, we have $u+pv\in[u,u+v]\subset\cI$ by convexity of $\cI$, and
thus $(u,u+pv)\in\cT_\cI$. Notice that $C^\Phi(\si_p)=p^2C^\Phi$.

\begin{lem}\label{le:cvx-4}
  Let $\Phi\in\cC_\cI$ and $p\in[0,1]$. Let $\si_p$ be as in and
  \eqref{eq:def-sigp}. The following inequalities hold on $\cT_\cI$, 
  \begin{equation}\label{eq:Ap-pA-Bp-pB}
    A^\Phi(\si_p)\leq pA^\Phi
    \text{\quad and\quad}
    B^\Phi(\si_p)\leq pB^\Phi,
  \end{equation}
  where $q:=1-p$. Moreover, $A^\Phi(\si_p)=p^2A^\Phi$ and
  $B^\Phi(\si_p)=p^2B^\Phi$ for \textbf{(P2)}. Let $\tau$ be as in
  \eqref{eq:def-tau}. Assume in addition that $\Phi''\in\cC_\cI$, then the
  following inequalities hold on $\cT_\cI$,
  \begin{eqnarray}
    pA^\Phi-A^\Phi(\si_p) 
    & \leq & pq\PAR{pA^\Phi(\tau)+qA^\Phi}; \label{eq:pA-Ap} \\
    A^\Phi(\si_p) 
    & \leq & \frac{1}{2}p^2q\,C^\Phi+p^3A^\Phi; \label{eq:Ap-C-A} \\
    A^\Phi(\tau(\si_p))) 
    & \leq & \frac{1}{2}p^2q\,C^\Phi+p^3A^\Phi(\tau); \label{eq:Atp-C-A} \\
    B^\Phi(\si_p) 
    & \leq & p^2q\,C^\Phi + p^3 B. \label{eq:Bp-C-B}
  \end{eqnarray}
\end{lem}

\begin{proof}
  The $A^\Phi$ part of \eqref{eq:Ap-pA-Bp-pB} is a rewriting of
  \eqref{eq:ent-twop}. For the $B^\Phi$ part of \eqref{eq:Ap-pA-Bp-pB}, we
  notice that $\Phi'$ is non-decreasing since $\Phi$ is convex, and thus
  $pv\Phi'(u+pv)\leq pv\Phi'(u+v)$ regardless of the sign of $v$, which gives
  the desired result. The inequality \eqref{eq:pA-Ap} is a rewriting of
  \eqref{eq:bern-mlsi-A}. Namely, $pA^\Phi(u,v)-A^\Phi(\si_p(u,v)) =
  \ENT{\cB(1,p)}{\Phi}{f}$ for any $(u,v)\in\cT_\cI$, where
  $(f(0),f(1)):=(u,u+v)$, whereas $\moy{\cB(1,p)}{A^\Phi(f,\rD f)} =
  pA^\Phi(\tau(u,v))+qA^\Phi(u,v)$. The inequalities \eqref{eq:Ap-C-A},
  \eqref{eq:Atp-C-A}, and \eqref{eq:Bp-C-B} follow from \eqref{eq:AB-int-C} by
  using the definition of $C^\Phi$ and a suitable Jensen inequality for
  $\Phi''$.
\end{proof}

\begin{rem}
  The bounds \eqref{eq:Ap-C-A} and \eqref{eq:Atp-C-A} become equalities for
  \textbf{(P2)}. However, \eqref{eq:Bp-C-B} is not sharp for \textbf{(P2)}.
  The bound $B^\Phi(\si_p)\leq pB^\Phi$ is optimal in the sense that for
  \textbf{(P1)}, we have $B^\Phi(u,pv)\sim pB^\Phi(u,v)$ at $v=+\infty$ for
  any $(p,u)\in(0,1)\times\cI$. However, $B^\Phi(\si_p)=p^2B^\Phi$ for
  \textbf{(P2)}; and in general
  $$
  \lim_{v\to0}v^{-2}B^\Phi(u,pv) %
  =p^2\lim_{v\to0}B^\Phi(u,v) %
  =p^2\Phi''(u).
  $$
  The same remark holds for $A^\Phi$ (up to a factor $2$ in the case
  \textbf{(P2)}).
\end{rem}

Some of the results of this section correct mistakes discovered by the author
in \cite{MR2081075} after publication. In contrary to what appears in
\cite[Page 330]{MR2081075}, \textbf{\cH(2')} does not imply \textbf{\cH(2)}.
Actually, \textbf{\cH(1)} and \textbf{\cH(2)} are equivalent and
\textbf{\cH(2')} should be removed from \cite{MR2081075}. In particular,
\cite[Remarks 8 and 11 pages 354-356]{MR2081075} should be replaced by Lemma
\eqref{le:cvx-1}. These corrections are minor and simplifying and do not
impact the results of \cite{MR2081075} at all.

\textbf{Acknowledgement.}
\emph{The 
  author is grateful to the anonymous referees of ESAIM/PS. Their remarks 
  and questions helped to improve the last version of the article. The author
  would also like to sincerely thank his friend and colleague Florent Malrieu
  for numerous stimulating discussions in Rennes, Saint-Malon sur Mel, and
  Toulouse.}

\providecommand{\bysame}{\leavevmode\hbox to3em{\hrulefill}\thinspace}
\providecommand{\MR}{\relax\ifhmode\unskip\space\fi MR }
\providecommand{\MRhref}[2]{%
  \href{http://www.ams.org/mathscinet-getitem?mr=#1}{#2} }
\providecommand{\href}[2]{#2}

\begin{center}
\hrule
\end{center}

Djalil \textsc{Chafa\"{\i}}. 
\textbf{Address.} UMR 181 INRA/ENVT Physiopathologie et
Toxicologie Experimentales, \'Ecole Nationale V\'et\'erinaire de Toulouse,
23 Chemin des Capelles, F-31076, Toulouse \textsc{Cedex} 3, France,
\textbf{and} UMR 5583 CNRS/UPS Laboratoire de Statistique et Probabilit\'es,
Institut de Math\'ematiques de Toulouse, Universit\'e Paul Sabatier, 118
route de Narbonne, F-31062, Toulouse, \textsc{Cedex} 4, France. 
\textbf{e-mail.} \url{mailto:chafai(AT)math.ups-tlse.fr.nospam}
\textbf{web-page.} \url{http://www.lsp.ups-tlse.fr/Chafai/}

\end{document}